\RequirePackage{fix-cm}
\documentclass[smallextended]{svjour3}
\usepackage{amsfonts,amsmath,amssymb,color}
\usepackage[round]{natbib} 
\smartqed  
\usepackage{graphicx}

\usepackage{color} %colori
\definecolor{bole}{rgb}{0.47, 0.27, 0.23}

\begin{document}

\title{Optimal control of the
mean field   equilibrium for a pedestrian tourists' flow model}
\titlerunning{Optimal control of the MFG equilibrium for a pedestrian
	tourists' flow model} % if too long for running head
\author{Fabio Bagagiolo \and Silvia Faggian \and Rosario Maggistro \and Raffaele Pesenti %etc.
}

\institute{F. Bagagiolo (corresponding author)\at
	Department of Mathematics, Universit\`{a} di Trento, Italy  \email{fabio.bagagiolo@unitn.it} \and 
	S. Faggian \at Department of Economy, Universit\`{a} Ca' Foscari Venezia, Italy \email{faggian@unive.it}   
	\and R. Maggistro \at
	 Department of Management, Universit\`{a} Ca' Foscari Venezia, Italy \email{rosario.maggistro@unive.it} \and R. Pesenti \at Department of Management, Universit\`{a} Ca' Foscari Venezia, Italy \email{pesenti@unive.it}
}

\date{Received: date / Accepted: date}
\maketitle

\begin{abstract}
Art heritage cities are popular tourist destinations but for many of them
overcrowding is becoming an issue. In this paper, we address the problem of
modeling and analytically studying the flow of tourists along the narrow
alleys of the historic center of a heritage city. We initially present a
mean field game model, where both continuous and switching decisional
variables are introduced to respectively describe the position of a tourist
and the point of interest that he/she may visit. We prove the existence of a
mean field   equilibrium. A mean field   equilibrium is Nash-type
equilibrium in the case of infinitely many players. Then, we study an
optimization problem for an external controller who aims to induce a
suitable mean field   equilibrium.
%\PACS{PACS code1 \and PACS code2 \and more}
\end{abstract}



\section{Introduction}

\label{intro} In the recent years, art heritage cities have experienced a
continuous growth of tourists' flow, to the point that overcrowding is
becoming an issue, and local authorities start implementing countermeasures.
In this paper we address the problem of modeling and analytically studying
the flow of tourists (or, more precisely, of daily pedestrian excursionists)
along the narrow alleys of the historic center of a heritage city. Starting
from the contents of \cite{bagpes}, we recast those
results into a mean field game with controlled dynamics on a network,
representing the attractive sites in the city and the possible paths to
reach them.

\begin{figure}[tbp]
\label{fig:circularnetwork} \centering
\includegraphics[scale=0.7]{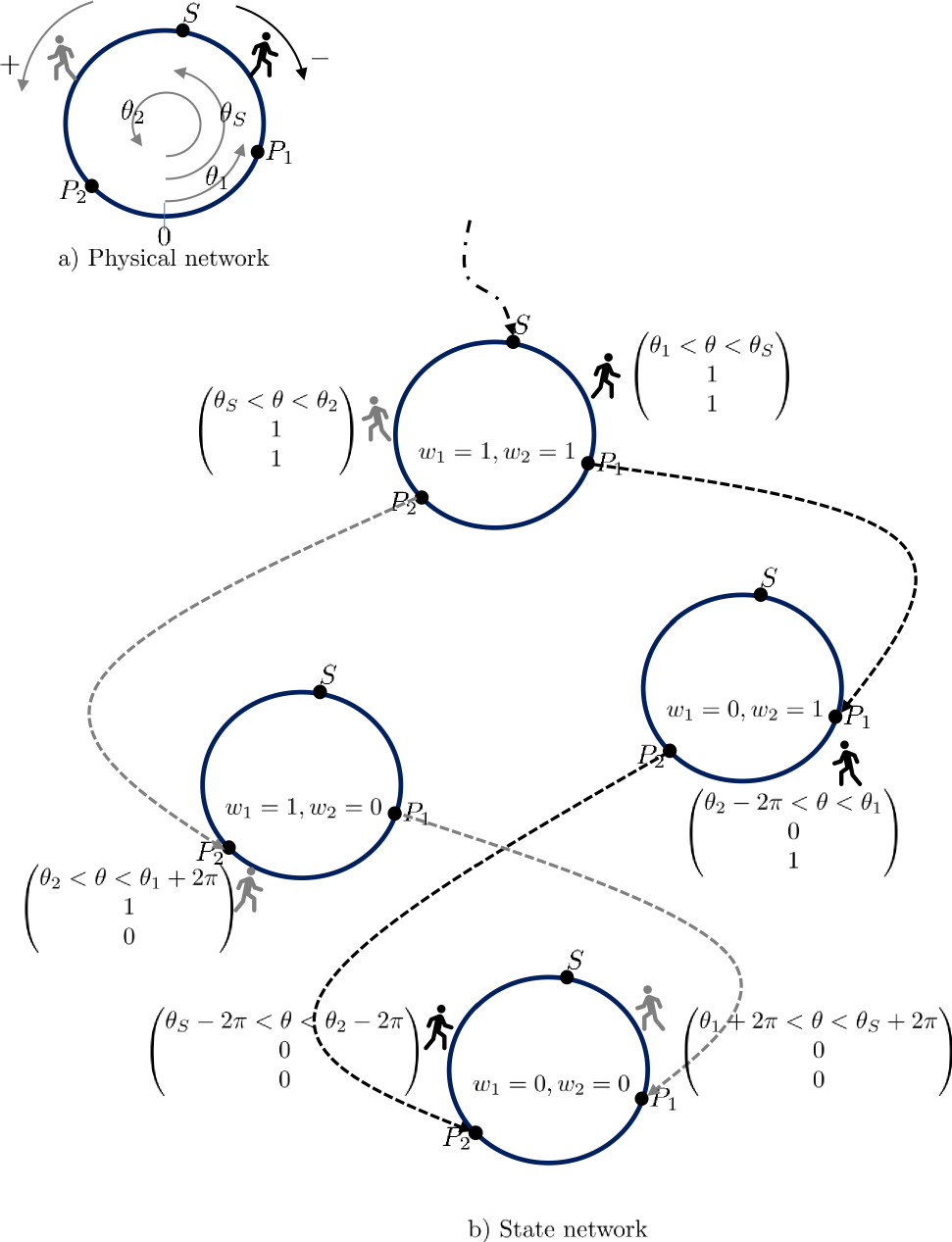}
\caption{a) the physical network of the paths inside the city,  {counterclockwise positively oriented}   (with the
three identified points: train station, $S$, and two attractions, $P_1, P_2$%
) and two tourists visiting the city following opposite directions. b) the
state network and the states $(\protect\theta, w_1, w_2)$ of the tourists
during their visits. Each singular circle-branch represents the network of
the city paths as seen by tourist with given values of the switching
variables $w_1, w_2$. The dashed arrows between the  four circles-branches
represent the switching of the four labels: at the beginning the label is $%
(1,1)$; when the attraction $P_1$ is reached (and hence visited), the label
switches to $(0,1)$, and similarly when $P_2$ is reached; the last is $(0,0)$
which holds when both attractions are reached. The point-dashed arrow
represents the external arrival flow in the station.}
\end{figure}

We here assume for simplicity that tourists have only two main attractions
to visit, as a context with more attractions can be equally treated. We
represent possible paths as a circular network containing three nodes: the
train station $S$ where tourists arrive in the morning and to which they
have to return in the evening; the first attraction $P_1$; the second
attraction $P_2$ (see Figure \ref{fig:circularnetwork}.a). The arrival flow
at the station is exhogenous, given by a continuous function $%
g:[0,T]\to[0,+\infty[$ representing, roughly speaking, the density of
arriving tourists per unit of time. The time by which everyone has to be
returned to the station after the tour is the final horizon $T>0$. 

A single tourist (agent),  {starting at position $\theta$ at time $t$}, controls its own dynamics, represented by the
equation 
\begin{equation}  \label{eq:thetau}
\begin{cases}
\theta^{\prime }(s)=u(s), &s\in ]{t},T]\\
\theta(t)=\theta &
\end{cases}
\end{equation}

\noindent where $\theta(s) \in\mathbb{R}$ is  {the} space-coordinate in the
network, and $s\mapsto u(s)\in \mathbb{R}$ is a measurable and locally integrable control, namely $u\in L^1_{loc}(t,T)$. {We represent the network by a circle (see Figure 1), and we
denote with $\theta _{S}$, $\theta _{1}$ and $\theta _{2}$ respectively the
position of the station, and of the attractions $P_{1}$ and $P_{2}$. In particular, we assume that 

\begin{equation}
\label{eq:thetas}
 {0<\theta_1<\theta_S<\theta_2<2\pi,}
\end{equation}

\noindent
which is not restrictive due to the circularity representation of the state space. Moreover (see Figure 1 and Figure 2), when convenient, we will identify $\theta$ and $\theta\pm 2\pi$.}

To each
tourist, we associate a time-varying label $(w_{1},w_{2})\in \{0,1\}\times
\{0,1\}$. For $i\in \{1,2\}$, $w_{i}(t)=1$ means that, at the time $t$, the
attraction $P_{i}$ has not been visited yet, and $w_{i}(t)=0$ that the
attraction has been already visited. The state of an agent is then
represented by the triple $(\theta ,w_{1},w_{2})$, where $\theta $ is a
time-continuous variable and $w_{1},w_{2}$ are switching variables. In the
following, we denote by $B=[0,2\pi ]\times \{0,1\}\times \{0,1\}$ the state
space of variables $(\theta ,w_{1},w_{2})$, and in particular we call 
\textit{(circle) branch} any subset $B_{{w}_{1},{w}_{2}}$ of $B$ which
includes the states $(\theta ,{{w}_{1},{w}_{2}})$, with $(w_{1},w_{2})$
fixed and $\theta $ varying in $[0,2\pi ].$ Such branches correspond to
edges of the switching networks in Figure 2, which represents the network of
Figure 1 in an equivalent way. 

While the evolution of $\theta$ is governed by~(\ref{eq:thetau}), switching
variables can only evolve from $1$ to $0$, that is 
\begin{equation*}
w_i = 
\begin{cases}
1 & \text{for} \ s \in [t,\tau_i] \\ 
0 & \text{for}\ s \in ]\tau_i, T],%
\end{cases}%
\end{equation*}
with $\tau_i\in[t,T]$ the first instant at which the agent reaches and
visits attraction $P_i$, $i \in \{1,2\}$. This evolution is represented by
the arrows and the labels in Figure~\ref{fig:circularnetwork}.b (see also
Figure~\ref{fig:switchingnetwork}). To exemplify, consider {in Figure~\ref%
{fig:circularnetwork}{.b} an agent that arrives at the station and visits $%
P_1$ first: its initial state is $(\theta_S,1,1) \in B_{1,1}$; at $\tau_1\in[%
0,T]$, $w_1$ switches from 1 to 0 so that, immediately after $\tau_1$, the
tourist's state $(\theta,0,1)$ belongs to the branch $B_{0,1}$} (see also
Figure 2).

%%%%% vecchia versione %%%%%%%%%
%consider a tourist that arrives in the station and  visits $P_1$ first, at time $\tau$. It has state $(\theta_S,1,1) \in B_{1,1}$ when it arrives at the station. Then, at time $\tau\in[t,T]$, its state $w_1$ switches, it was equal to 1 in $[t,\tau_1]$, it assumes value $w_1=0$ in a left-open interval $]\tau_1,T]$. Immediately after  $\tau_1$, the tourist's state switches to the branch  $B_{0,1}$.} See also Figure 2.
%%%%%%%

The cost to be minimized by every agent 
%when it starts its tour from the station at time $t$, 
takes into account: i) the hassle of running during the tour; ii) the pain
of being entrapped in highly congested paths; iii) the frustration of not
being able to visit some attractions; iv) the disappointment of not being
able to reach the station by the final horizon $T$. Such a cost can be
analytically represented by

\begin{equation}  \label{eq:cost}
\begin{array}{ll}
J(u; t,\theta, w_1,w_2)=\int_t^T\left(\frac{u(s)^2}{2}+\mathcal{F}^{w_1(s),w_2(s)}(\mathcal{M}%
(s))\right)ds+ &  \\ 
\ \ \ \ \ \ \ \ \ \ \ \ \ \ \ \ \ \ \ \ \ \
+c_1w_1(T)+c_2w_2(T)+c_S\xi_{\theta=\theta_S}(T) & 
\end{array}%
\end{equation}

Here, $c_1,c_2,c_S>0$ are fixed, and $\xi_{\theta=\theta_S}(s)\in\{0,1\}$
and it is equal to $0$ if and only if $\theta(s)=\theta_S$. In (\ref{eq:cost}%
), the quadratic term inside the integral stands for cost i) while the other
term stands for the congestion cost ii); the following two addenda stand for
costs iii); the last addendum stands for cost iv). In particular, the
congestion cost $\mathcal{F}^{w_1(s),w_2(s)}(\mathcal{M}(s))$ is
instantaneously paid by an agent whose switching label at time $s$ is $%
(w_1(s),w_2(s))$, being $\mathcal{M}(s)$ the actual distribution of the
agents. For any $(w_1,w_2)\in\{0,1\}\times\{0,1\}$, $\mathcal{F}^{w_1,w_2}$
is a positive function defined on the set of all admissible distribution of
agents.
 
 The problem here treated differs from that  in  \cite{bagpes} in the respect that  the discontinuous final cost $c_S\xi_{\theta=\theta_S}(T)$  in  
(\ref{eq:cost})  is replaced by  
a smooth cost  $c_3(\theta-\theta_S)^2$. As done in \cite{bagpes}, the problem is considered  as a mean field game, and
the existence of a mean field  equilibrium
  is proven under suitable
assumptions. A mean field   equilibrium is a time-varying distribution of
agents on the network, $s\mapsto\mathcal{M}^*(s)$ for $s\in[0,T]$
generating, when plugged in (\ref{eq:cost}), an optimal control $%
u^*$  and associated optimal
trajectory, for any agent starting at the station at time $t\in[0,T]$,
yielding again the time-varying distribution $\mathcal{M}^*$. A mean field
  equilibrium can be seen as a fixed point, over a suitable set of
time-varying distributions, of a map of the form 
\begin{equation}  \label{eq:fixedpoint}
\mathcal{M}\longrightarrow u_{\mathcal{M}}\longrightarrow \mathcal{M}_{u_{%
\mathcal{M}}}
\end{equation}
where $u_{\mathcal{M}}$ is the optimal control when $\mathcal{M}$ is
inserted in (\ref{eq:cost}), and $\mathcal{M}_{u_{\mathcal{M}}}$ is the
corresponding evolution of the agents' distribution when all of them are
moving implementing $u_{\mathcal{M}}$ as control. Of course, the problem
must be coupled with an initial condition for the distribution $\mathcal{M}$%
, whereas the boundary condition is represented by the incoming flow $g$.

We remark that the concept of mean field equilibrium is of Nash type.
Indeed, in the case of a large number of agents (even infinitely many, as in
the case of a mean field game) every single agent is irrelevant, the single
agent has measure zero, it is lost in the crowd. Hence, at equilibrium, for
a single agent is not convenient to unilaterally change behavior, because
such a single choice will not change the mean field $\mathcal{M}$, and the
agent will not optimize his behavior.

That said, the goal of the present paper is twofold. First we amend some
stringent assumptions that were made in \cite{bagpes} to prove the existence
of a mean field   equilibrium. Then, we study a possible optimization
problem for an external controller who aims to induce a suitable mean field
 equilibrium. 
We suppose that the external controller (the city
administration, for example) may act on the congestion functions $\mathcal{F}%
^{w_1,w_2}$, choosing them among a suitable set of admissible functions.

We also note that the presence of a discontinuous final cost is by-passed using a  characterization of the value function at some significant states, identifying at the same time a finite set of optimal controls.

Consider as an example the historical center of the city of Venice, Italy.
Tourists typically enter the city at the train station and from there they
have two main routes, a shorter and a longer one, to major monuments. In
recent years, the shorter route has been particularly congested on peak
days. For this reason, local authorities have introduced both some gates to
slow down the access to the shorter route and some street signs to divert at
least part of the tourist flow along the longer route.

Mean field games theory goes back to the seminal work by  \cite{LLions} (see also \cite{HCAINMAL}), where the new
standard terminology of mean field games was introduced. This theory
includes methods and techniques to study differential games with a large
population of rational players and it is based on the assumption that the
population influences the individuals' strategies through mean field
parameters. Several application domains such as economic, physics, biology
and network engineering accommodate mean field game theoretical models
 \citep{AcCamDolc,LachSalTur}. In particular, models to the
study of dynamics on networks and/or pedestrian movement are proposed, for
example, by \cite{CCMar}, 
by \cite{CPTos},  by \cite{CDeMTos}, and by
 \cite{BBMZop}. 

Very recently  \cite{LLions2018} have introduced a new class of mean field games
to model situations involving one dominant agent and a large group of small players.
    Our work can be framed in this new class when it deals with the possible role of local authorities in the optimization of tourist flows.

\smallskip

The problem here treated is also a routing problem 
(for a particular case
of pedestrian movement with possible multiple destinations see, for example, \cite{Hoogendoorn}).
 Note that different
strategies were proposed to control the roadway congestion, such as variable
speed limits \citep{HegyiSchutter1}, ramp
metering \citep{GomesHoro} or signal control \citep{BrianYun}. However, such mechanisms neither consider the
agents' perspective, nor affect the total amount of vehicles/people. A
significant research effort was done to understand the agents' answer to
external communications from intelligent traveler information devices \citep{Khattak,Mahamassani} and, in particular, to study the effect of such technologies
on the agents' choice and on the dynamical properties of the transportation
network \citep{ComoSavla}. Moreover
it is well known that when individual agents make their routing decisions to
minimize their own experienced delays, then the overall network congestion
is considerably higher than the congestion resulting from a central planner
directing traffic explicitly \citep{Pigou}. From that the idea to include in our problem
an external controller inducing a suitable mean field   equilibrium.

 Note also that in the specific case of vehicular congestion,
tolls payment is considered to influence drivers to make routing choices
that result in globally optimal routing, namely to force the Wardrop
equilibrium to align with the system optimum network flow 
\citep{Smith,Morrison,Dial,Cole}. The Wardop equilibrium, proposed by 
\cite{Wardrop}, is a configuration in which the perceived cost associated to
any source-destination path chosen by a nonzero fraction of drivers does not
exceed the perceived cost associated to any other path. This is
a stationary  concept, whose
continuous counterpart was recently
developed by \cite{Carlier1} and by \cite{CarlierSant}. 

Not only does the latter fit the situation of pedestrian congestion 
%(that in this paper we indeed treat using a mean field model) 
but also is useful to analyse large-scale traffic problems, in which one
only wants to identify the average value of the traffic congestion in
different zones of a large area. %
% It is described through the formalism of measure on the set of paths, which is a classical tool in transportation theory. 
Actually, also the models using that continuous framework are essentially stationary, as 
% Moreover in this continuous framework some general minimization problems, whose optimality conditions are expressed by a partial differential equation (PDE), can be consider. Notice anyway that the equations involved in these formulations are elliptic PDE's, and no variable playing the role of time appears in them. 
% This because the model is essentially stationary 
% (time is only fictitious,
%and traffic intensity is considered as a function of the state and not of the time): 
they only account for a sort of cyclical movement, where every
path is constantly occupied by the same density of vehicles, since those
reaching the destination are immediately replaced by others. This is the
essential difference with respect to the mean field models. Indeed in the
latter, due to the explicit presence of time, the optimal evolution is given
by a system coupling a transport equation and an Hamilton-Jacobi equation. %

%%%%%%%%%%%%%%%%%%%%%%%%%%%%%%%%%%%%%%%%%%%%%%%%%
%%%%%%%%%%%%%%%%%%%%%%%%%%%%%%%%%%%%%%%%%%%%%%%%%
%%%%%%%%%%%%%%%%%%%%%%%%%%%%%%%%%%%%%%%%%%%%%%%%%
%%%%%%%%%%%%%%%%%%%%%%%%%%%%%%%%%%%%%%%%%%%%%%%%%

\section{Preliminary results}

\label{sec:1} 
 
As anticipated in the introduction, the problem here treated differs from 
that in  \cite{bagpes}, under the respect that the final cost $c_S\xi_{\theta=\theta_S}(T)$  in  
(\ref{eq:cost}) replaces the smooth cost $c_3(\theta-\theta_S)^2$ of  \cite{bagpes}. 
Hence, we here recall some of the results from the previous paper
that still hold true in our case, and present some new ways to approach
the problem.

 We define the value
function $V(\theta,t, w_1,w_2)$   as the infimum over the measurable
controls of the cost (\ref{eq:cost}) faced by an agent starting its
evolution at the state $(\theta,w_1,w_2)$, at time $t$, that is

\[
V(\theta,t,w_1,w_2)=
\inf_{u\in L^1_{loc}(t,T)}J(u; t,\theta, w_1,w_2)
\]

Note that when the distributional evolution $t\mapsto%
\mathcal{M}(t)$, i.e. the density of the agents, is initially given, then  $V$ does not depend on $\mathcal{M}$.  
As well as in \cite{bagpes}, the following two sets of equations are considered:
\begin{itemize}
\item[(a)] the set of  Hamilton-Jacobi-Bellman equations (one for each branch $B_{w_1,w_2}$) associated to the problem, 
and fully described through (\ref{eq:HBJXX}) in Appendix B (corresponding to equations
 {(3)--(6)} in 
\cite{bagpes});
\item[(b)] the associated set of transport equations (one for each branch), fully described through (\ref{eq:transportXX}) in Appendix B (corresponding to  equations  {(8a)--(8d)} in \cite{bagpes}).
\end{itemize}

\begin{figure}[tbp]
\label{fig:switchingnetwork} %
\includegraphics[scale=0.45]{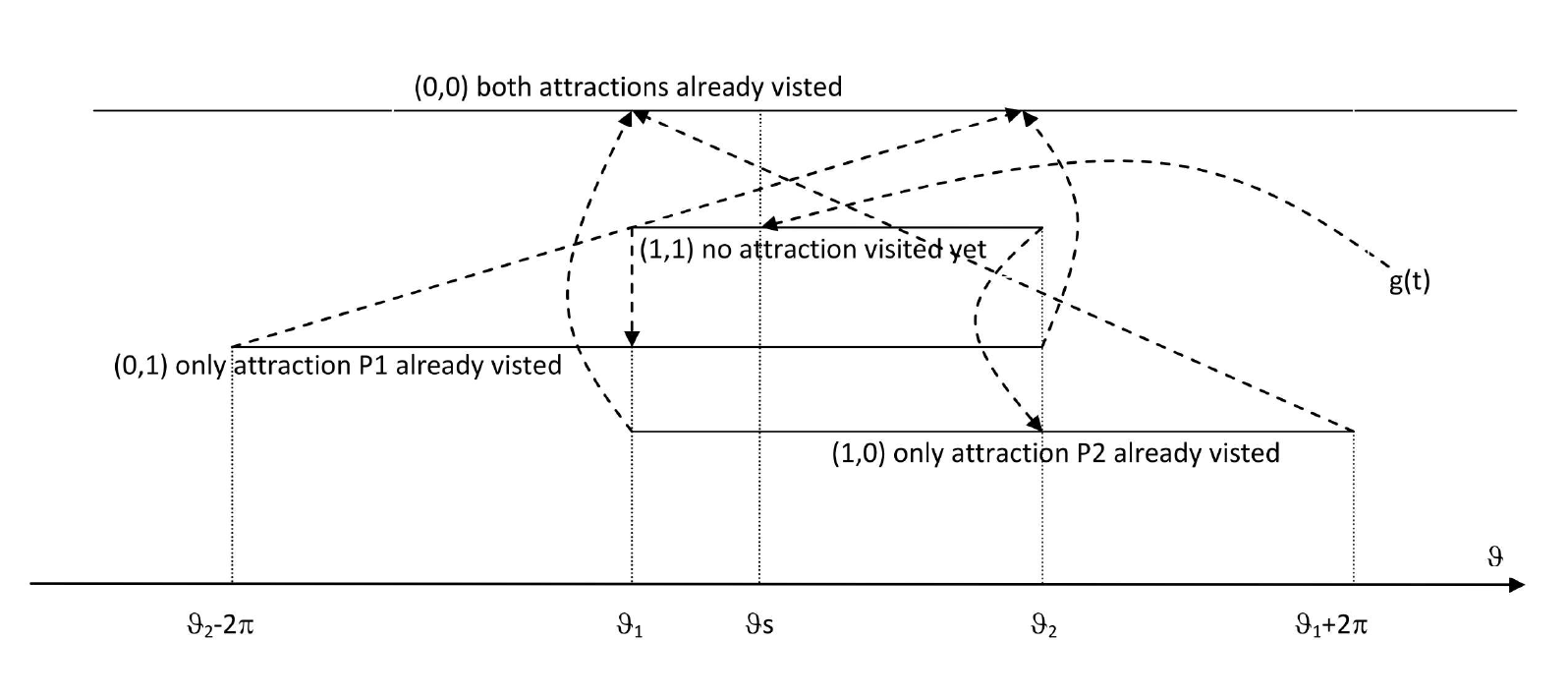}
\caption{Switching cross branches representation}
\end{figure}

Regarding (a), note that in (\ref{eq:HBJXX}) boundary conditions are given by the value, at  the switching points, of the solution in the consecutive branch. 
This fact comes from the dynamic programming principle, and from interpreting the
optimal control problem on a single branch as an exit-time problem \citep{bagdan}. More in detail (see Figure \ref{fig:switchingnetwork}), the exit cost from $B_{1,1}$ at time $\tau$ is
given by $V(\theta_1,\tau,0,1)$ at $\theta_1$, and by $V(\theta_2,\tau,1,0)$
at $\theta_2$. The exit costs from $B_{0,1}$ and from $B_{1,0}$ are,
respectively, $V(\theta_2,\tau,0,0)$ and $V(\theta_1,\tau,0,0)$ (on both
exit points of the branch). On the final branch $B_{0,0}$ the problem is
reduced to reaching the station at the least cost within the final time $T$.
\medskip

Note also that, when  the distribution $\mathcal{M}$ is given, the optimal
control problem faced by an agent on a branch 
is a rather standard combination of a reaching target (the station) problem
and suitable exit time problems. In particular, the optimal control problem
on $B_{0,0}$ can be solved, and $V(\cdot,\cdot;0,0)$ calculated,
independently from the evolution of the system on the other branches. Once
the problem on $B_{0,0}$ is solved, one may proceed backwards and solve the
problems on $B_{0,1}$ and $B_{1,0}$ and, eventually, on $B_{1,1}$.

%Once the control problem is solved on every branch, the optimal transport problem  described by \cite[(1.8a)--(1.8d)]{bagpes} can be tackled. 

Regarding (b), the {transport equations}  (\ref{eq:transportXX}) are solved,  for every $(w_1,w_2)\in\{0,1\}^2$ and for every $t\in[%
0,T]$, by  a  time dependent measure $m^{w_1,w_2}(\cdot,t):[0,2\pi]\to[0,+\infty[$  representing the agents' distribution on the branch $B_{w_1, w_2}$
at time $t$. It is assumed that the city network is initially empty of
tourist, that is, $m^{w_1,w_2}(\cdot,0) = 0$ for $(w_1,w_2)\in\{0,1\}^2$. By conservation of mass principle, $\mathcal{M}$ $%
=(m^{1,1},m^{0,1},m^{1,0},m^{0,0})$ satisfies 
\begin{align*}
\int_Bd\mathcal{M}(t)=&\int_{B_{1,1}}dm^{1,1}(t)+\int_{B_{1,0}}dm^{1,0}(t) +
\\
&+\int_{B_{0,1}}dm^{0,1}(t) +\int_{B_{0,0}}dm^{0,0}(t)=\int_0^tg(s)ds, \quad
t \in [0,T].
\end{align*}
The connection between the solutions $m^{w_1,w_2}$ in
the different branches is obtained again through boundary conditions in  (\ref{eq:transportXX}).

\bigskip

As in  \cite{bagpes}, we here assume the following hypotheses:
\begin{itemize}
\item[(H1)]  $g:[0,T]\to[0,+\infty[$ is a Lipschitz continuous function;
\item[(H2)] $m^{w_1,w_2}$ are  continuous functions of time into the set ${\cal B}_{w_1,w_2}$ of Borel measures on the corresponding branch $B_{w_1,w_2}$, endowed with the weak-star topology, and ${\mathcal M}(0)=0$;
\item[(H3)] $t\mapsto \mathcal{F}^{w_1,w_2}(\mathcal{M}(t))$   continuous and
bounded for all  $(w_1,w_2)\in\{0,1\}^2$; in particular, ${\cal F}^{w_1,w_2}$ is continuous as function from 
${\cal B}_{1,1}\times\dots\times{\cal B}_{0,0}$ to $\mathbb R$;
\item[(H4)]  $\mathcal{F}^{w_1,w_2}$ does not depend explicitly on state variable $\theta$.
\end{itemize}

Note that ${\mathcal M}(0)=0$ in  (H2) means that no one is around the city at $t=0$, while (H4) means that all agents in the same branch at the same instant equally suffer
the same congestion pain.

We recall that in \citet[Theorem 1]{bagpes} it is proven, under the previous assumptions, that the value function $V$ coincides in every branch with the unique (continuous) viscosity solutions of the system of Hamilton-Jacobi-Bellman equations (\ref{eq:HBJXX}).
Moreover, for every branch $B_{w_1,w_2}$,   the optimal
feedback control $u^*$ needs to satisfy for all $t$ and $\theta$
\begin{equation}  \label{eq:optimalfeedback}
u^*(\theta,t,w_1,w_2)=-V_\theta(\theta,t,w_1,w_2).\end{equation}

\begin{remark}\label{ccp}
Note that (H1) (H3) (H4) also imply that the control choice made by agents at states $(\theta_S,1,1)$, $(\theta_1,0,1)$, $(\theta_2,1,0)$, $(\theta_1,0,0)$ and $(\theta_2,0,0)$ (that from now on we call \textit{significant states}) does not change as long as the agent remains in the same branch, and is constant in time.
  To see this, we observe beforehand that if an agent moves from $\theta'$ at time $t'$ and reaches $\theta''$ at time $t''>t'$, its cost $\frac 1 2 \int_{t'}^{t''}u^2(\tau)d\tau$ is minimized, by Jensen's inequality and \eqref{eq:thetau}, when the control is chosen constant and equal to  the mean speed $ u^*(\tau)\equiv \frac{\theta''-\theta'}{t''-t'}$. This fact excludes the possibility  - recall that ${\cal F}^{1,1}$ does not explicitly depend on the space variable $\theta$ -  that an  optimally behaving agent persists at a state (i.e, chooses $u=0$) along a positive time interval and  moves later, as this behaviour would be worse than choosing the mean speed from the beginning. 
As a consequence of the above argument, an agent has the choice,  at every significant state, of either stay still  forever   or move to the next significant state at  constant  speed. 
 Such a characterization of the optimal controls is used in the next subsections.
\qed 
\end{remark}
%%%%%%%%%%%%%%%%%%%%%%%%%%%%%%%%%%%%
%%%%%%%%%%%%%%%%%%%%%%%%%%%%%%%%%%%%
%%%%%%%%%%%%%%%%%%%%%%%%%%%%%%%%%%%%
%%%%%%%%%%%%%%%%%%%%%%%%%%%%%%%%%%%%

\subsection{Value Functions}
 The presence of the discontinuous final cost $c_S\xi_{\theta=\theta_S}(T)$ in (\ref{eq:cost}) calls for a characterization of the value function  different than  (\ref{eq:HBJXX}), and contained into  \eqref{eq:V00}  \eqref{eq:V01}  \eqref{eq:V10}  \eqref{eq:V11}. We remark that, along the process, we  will also easily detect all  available optimal controls in view of Remark \ref{ccp}.

An agent standing at $(\theta_i,0,0)$ at time $t\in[0,T]$, with $i\in\{1,2\}$,  has two possible
choices (meaning, the optimal behavior may only be one of the following
two): either staying at $\theta_i$ indefinitely or moving to reach $\theta_S$
exactly at time $T$ (it is not optimal to reach $\theta_S$ before $T$ and
wait there for a positive time length, as pointed out in \cite{bagpes}). The  controls among which the agent chooses are then, respectively
\begin{equation}\label{oc1}
{u^{0,0}_0(t)}\equiv0,\ \ u^{0,0}_1(t)=\pm\frac{\theta_S-\theta_1}{T-t},\ \ u^{0,0}_2(t)=\pm\frac{\theta_S-\theta_2}{T-t}.
\end{equation} Here $\theta_S-\theta_i$ stands
for the length of the minimal path between $\theta_S$ and $\theta_i$ in the
circular network, while the sign is chosen consistently with (\ref{eq:thetau}%
) in running that path. In particular, they are constant controls, determined by the arrival time ($T$ in this case) and the distance to be run.
Hence, given the cost functional (\ref{eq:cost}), we derive 
\begin{equation}  \label{eq:V00}
 V(\theta_i,t,0,0)=\min\left\{c_S,\frac{1}{2}\frac{%
(\theta_S-\theta_i)^2}{T-t}\right\}+\int_t^T\mathcal{F}^{0,0}\, ds
\end{equation}
(we do not display the argument of $\mathcal{F}^{0,0}$, that being the fixed
distribution $\mathcal{M}$). 
Instead an agent standing at $(\theta _{1},0,1)$ at time $t$ has three
possible choices: staying at $\theta _{1}$, moving to reach $\theta _{S}$ at
time $T$; moving to reach $\theta _{2}$ at $\tau \in ]t,T]$. Accordingly, the possible choices for a control are, respectively
\begin{equation}\label{oc2}
{u^{0,1}_0(t)}\equiv0,\ \ u^{0,1}_1(t)=  \pm \frac{\theta _{S}-\theta _{1}}{T-t},\ \ u^{0,1}_2(t)= \pm \frac{\theta _{2}-\theta _{1}}{\tau -t}.
\end{equation} 
 Hence, 
\begin{align}\label{eq:V01}
\displaystyle
V(\theta _{1},t,0,1)&=\min \left\{c_{2}+c_{S}+\int_{t}^{T}\mathcal{F}^{0,1}\,ds, c_{2}+\frac{1}{2}\frac{(\theta _{S}-\theta _{1})^{2}%
}{T-t}+\int_{t}^{T}\mathcal{F}^{0,1}\,ds,\right. \nonumber\\
& \left.\inf_{\tau\in ]t,T]}  \left\{\frac{1}{2}\frac{(\theta _{2}-\theta _{1})^{2}}{\tau
-t}+\int_{t}^{\tau }\mathcal{F}^{0,1}\,ds+V(\theta _{2},\tau ,0,0)\right\}\right\}   
\end{align}
% \color{red}
%and, on the whole branch $B_{0,1}$ $$V(\theta,t,0,1)=.....?????$$ , where\color{black}

Similarly at  $(\theta _{2},1,0)$ the control is chosen among
\begin{equation}\label{oc3}
{u^{1,0}_0(t)}\equiv0,\ \ u^{1,0}_1(t)=  \pm \frac{\theta _{S}-\theta _{2}}{T-t},\ \ u^{1,0}_2(t)= \pm \frac{\theta _{2}-\theta _{1}}{\tau -t}.
\end{equation} yielding
\begin{align}\label{eq:V10}
\displaystyle
{V(\theta _{2},t,1,0)}&=\min \left\{c_{1}+c_{S}+\int_{t}^{T}\mathcal{F}^{0,1}\,ds, c_{1}+\frac{1}{2}\frac{(\theta _{S}-\theta _{2})^{2}%
}{T-t}+\int_{t}^{T}\mathcal{F}^{0,1}\,ds,\right. \nonumber\\ 
& \left.\inf_{\tau\in ]t,T]}  \left\{\frac{1}{2}\frac{(\theta _{2}-\theta _{1})^{2}}{\tau
-t}+\int_{t}^{\tau }\mathcal{F}^{0,1}\,ds+V(\theta _{1},\tau ,0,0)\right\}\right\}. 
\end{align}
% and \color{red} $V(\theta,t,1,0)\equiv .................$ on the whole branch $B_{1,0}$. \color{black}

Finally, an agent standing at $(\theta _{S},1,1)$ at time $t$ may choose: to
stay there, to reach $\theta _{1}$ at a certain $\tau \in ( t,T]$, to
reach $\theta _{2}$ at a certain $\eta \in ( t,T]$. 
{At  $(\theta _{S},1,1)$ the control is chosen among
\begin{equation}\label{oc4}
u^{1,1}_0(t)\equiv0,\ \ u^{1,1}_1(t)=  \pm \frac{\theta _{S}-\theta _{1}}{\tau-t},\ \ u^{1,1}_2(t)= \pm \frac{\theta _{S}-\theta _{2}}{\eta -t}.
\end{equation}}
Consistently, one
has
\begin{align}\label{eq:V11}
\displaystyle
V&(\theta _{S},t,1,1)=\min \left\{ c_{1}+c_{2}+\int_{t}^{T}\mathcal{F}^{1,1}ds,         \inf_{\tau\in]t,T]}   \frac{1}{2}\frac{(\theta _{1}-\theta _{S})^{2}}{\tau -t}%
+\int_{t}^{\tau }\mathcal{F}^{1,1}\,ds   +             \right. \nonumber\\ 
& \left.+V(\theta _{1},\tau ,0,1), \inf_{\eta\in ]t,T]}\frac{1}{2}\frac{(\theta _{2}-\theta _{S})^{2}}{\eta
-t}+\int_{t}^{\eta }\mathcal{F}^{1,1}\,ds+V(\theta _{2},\eta ,1,0)    \right\}. 
\end{align}

Note that, the  optimal controls described in (\ref{oc1}), (\ref{oc2}), (\ref{oc3}), (\ref{oc4})
are detected  along with the  arrival times $\tau $ and $\eta $ along the minimization process carried on in  (\ref{eq:V00}), (\ref{eq:V01}), (\ref{eq:V10}), (\ref{eq:V11}).
Note again that, when $\cal M$ is given, the construction of the optimal controls is  performed  backwardly, starting from the minimization problem (\ref{eq:V00}).
\medskip

We summarize the previous discussion as follows.

%%%%%%%%%%%%%%%%%%%%%%%%%%%%%%%%%%%%%%%%%%%%%%%%%

	\begin{theorem} \label{th:valuefun} Assume that $\cal M$ is given and that (H1)--(H4) hold. Then the value function $V : B \times [0,T] \to \mathbb{R}$, at the significant states $(\hat \theta,w_1,w_2)$, with $\hat\theta\in\{\theta_S,\theta_1,\theta_2\}$, and for all $t\in[0,T]$,
	is determined through (\ref{eq:V00})-(\ref{eq:V11}). 
In addition, $V$% \deleted{is the unique viscosity solution to the Hamilton-Jacobi-Bellman equations (\ref{eq:HBJXX}) and it}
is Lipschitz continuous with respect to time, and its Lipschitz constant is independent of $\mathcal M$.
\end{theorem}

\begin{proof} The proof that the value function is characterized  by  (\ref{eq:V00})-(\ref{eq:V11}) is contained in the previous discussion. 
	 It is sufficient to  prove the Lipschitz continuity with respect to time at the significant states,  as implied by (H4), and by the fact that once a control is chosen at one of those states, then it is  maintained up to the   exit from the branch. 
  Note now that  (H3) implies there exists a positive constant $k_1$, independent of $\mathcal M$, such that
\begin{equation*}
\left\vert \int_{t}^{\tau }{\mathcal F}^{w_{1},w_{2}}({\mathcal M(s)})ds\right\vert \leq
\left\Vert {\mathcal F}^{w_{1},w_{2}}({\mathcal M(\cdot )})\right\Vert _{\infty }\left\vert \tau
-t\right\vert \le k_{1}\left\vert \tau -t\right\vert \leq k_{1}T
\end{equation*}
for all $%
\tau \in \lbrack t,T],$ and for all $(w_{1},w_{2})$.
That in particular implies that $V(\theta_i,t,0,0)$ given by (\ref{eq:V00})  is Lipschitz
continuous in $t$ and bounded in absolute value by a constant $k_{2}$, also not depending on ${\mathcal M}$.
Proceding backwards, this implies that also $V(\theta_1,t,0,1)$ given by (\ref{eq:V01}) is
Lipschitzian, as the infimum with respect to $\tau$ of appearing in (\ref{eq:V01})  can be
computed on $[t+h,T]$, with $h$ such that%
\begin{equation}
\label{h}
\frac{1}{2}\frac{(2\pi )^{2}}{h}+k_{2}\leq c_{1}+c_{2}+c_{S}.
\end{equation}
Arguing similarly, one proves the Lipschitz continuity on the branches $B_{1,0}$ and $B_{1,1}.$ \qed
\end{proof}

\subsection{Masses and Flows of Agents}
Along with significant states, we consider the arrival flows at
such states (see Figure \ref{fig:switchingnetwork}): the given external arrival {flow} at the station $g$, and the four  flows  $g_{0,1}$ at $%
(\theta_1,0,1)$, $g_{1,0}$ at $(\theta_2,1,0)$; $g_{1,2}$ at $(\theta_2,0,0)$,
and $g_{2,1}$ at $(\theta_1,0,0)$. The flows need satisfy the
following conservation constraints for all $t\in[0,T]:$ 
\begin{equation}  \label{eq:conscond}
\begin{array}{ll}
\displaystyle \int_0^tg(\tau)d\tau\ge\int_0^tg_{0,1}(\tau)d\tau+%
\int_0^tg_{1,0}(\tau)d\tau, &  \\ 
\displaystyle \int_0^tg_{0,1}(\tau)d\tau\ge\int_0^tg_{1,2}(\tau)d\tau, &  \\ 
\displaystyle \int_0^tg_{1,0}(\tau)d\tau\ge\int_0^tg_{2,1}(\tau)d\tau. & 
\end{array}%
\end{equation}
 Note that the flow functions $g_{0,1},g_{1,0},g_{1,2}$ and $g_{2,1}$, as well as the given flow $g$, are time densities entering the switching states. They generate  spatial densities along the branches $B_{w_1,w_2}$, and spatial densities transform once again into time densities at the subsequent switching point. The reader may find the thorough description of spatial and temporal components of flow functions in  Appendix {A}.
\medskip

Denoting by $\rho^{w_1,w_2}(t)$ the actual total mass of agents on the
branch $B_{w_1,w_2}$, we then have 
\begin{equation}  \label{eq:totalmass}
\begin{array}{ll}
\displaystyle \rho^{1,1}(t)=\int_0^tg(\tau)d\tau-\int_0^tg_{0,1}(\tau)d\tau-%
\int_0^tg_{1,0}(\tau)d\tau, &  \\ 
\displaystyle \rho^{0,1}(t)=\int_0^tg_{0,1}(\tau)d\tau-\int_0^tg_{1,2}(\tau)d%
\tau, &  \\ 
\displaystyle \rho^{1,0}(t)=\int_0^tg_{1,0}(\tau)d\tau-\int_0^tg_{2,1}(\tau)d%
\tau, &  \\ 
\displaystyle \rho^{0,0}(t)=\int_0^tg_{1,2}(\tau)d\tau+\int_0^tg_{2,1}(\tau)d%
\tau. & 
\end{array}%
\end{equation}
and, with the notation of the previous section, 
\begin{equation}\label{rho}
\rho^{w_1,w_2}(t)=\int_{B_{w_1,w_2}}{dm^{w_1,w_2}(t)}.
\end{equation}
We denote by $\rho=(\rho^{1,1},\rho^{0,1},\rho^{1,0},\rho^{0,0})$ the vector of masses in
different branches.

\begin{remark}\label{rholip}
%Note that the flow functions $g_{i,j}$ and $\rho^{w_1,w_2}$ are bounded: this is straightforward from (\ref{eq:conscond}),(\ref{eq:totalmass}).  Moreover 
%\begin{equation}
%\label{K}0\le \rho^{w_1,w_2}(t)\le K\equiv\int_0^Tg(s)ds.
%\end{equation} 
%The  functions $\rho^{w_1,w_2}$ are also all Lipschitz continuous, with the same Lipschitz constant $\tilde L$, as $g$ and $g_{i,j}$ are uniformly bounded.\qed

All the possible functions $\rho^{w_1,w_2}$ are obviously uniformly bounded by

\begin{equation}
\label{K}
0\le \rho^{w_1,w_2}(t)\le K\equiv\int_0^Tg(s)ds.
\end{equation}

\noindent

Moreover, since the optimal controls are necessarily equi-bounded by a constant depending only on the parameters of the problem 
(because of the presence of the term $u^2/2$ in the cost), and since the flow $g$, entering at $\theta_S$ in the branch $B_{1,1}$,  is continuous and hence bounded, by results on the transport equations (see for example the unpublished notes by Cardaliaguet (2013):
https://www.ceremade.dauphine.fr/$\scriptstyle\sim$cardaliaguet/MFG20130420.pdf), we get that $\rho^{1,1}$ is Lipschitz continuous, with Lipschitz constant depending only on the parameters of the problem. Hence it has bounded derivatives almost everywhere. Differentiating with respect to $t$ the first line of (\ref{eq:totalmass}), we obtain that also the flow functions $g_{0,1}$ and $g_{1,0}$ (which are positive) are almost everywhere bounded, by a constant only depending on the parameters of the problem. Proceeding in this way in the subsequent branches, we obtain the uniform boundedness of all flow functions and the uniform Lipschitz continuity of all possible functions $\rho^{w_1,w_2}$ in the other branches.  {Let us denote by $\tilde L$ the common Lipschitz constant}. We also point out that (as explained in  Appendix {A}), the flow functions $g_{i,j}$ are necessarily quite regular, because moving agents (i.e. agents not using null control) cannot accumulate at any time less than the final time $T$. Moreover, they have bounded integrals by (\ref{eq:conscond}).

\end{remark}

From now on, we assume that $\mathcal{F}^{w_1,w_2}$ only depends on $\rho$, instead of $\cal M$, and in particular on $\rho^{w_1,w_2}$ only. 
With an abuse of notation, we keep indicating the function by  $\mathcal{F}^{w_1,w_2}$ and replace  (H3) with the following (simplifying) assumption.
\begin{itemize}
\item[(H$3^\prime$)] $\mathcal{F}^{w_1,w_2}:[0,{K}]\to[0,+\infty[$  {are Lipschitz} continuous, and \\
$\mathcal{F}^{w_1,w_2}({\cal M}(t))=\mathcal{F}^{w_1,w_2}(\rho^{w_1,w_2}(t)),\, \forall t\in[0,T], \forall (w_1,w_2)\in\{0,1\}^2.$
\end{itemize}
%Indeed, (H$2^\prime$) and Remark \ref{rholip} imply (H2). 
 Assumption (H$3^\prime$) implies that  congestion costs  $\mathcal{F}^{w_1,w_2}$
	only depend on the total mass $\rho^{w_1,w_2}$ of agents on the single branch~$B_{w_1,w_2}$, 
	instead that on agents' distribution~$\mathcal{M}$ or on the entire vector $\rho$. This simplification, although less realistic from
	a modeling perspective, allows to completely treat the model. 
	Note that flows from and to the station $S$
	are, in the practice, distributed in different times, hence assumption (H$3^\prime$), although simplistic, is still
	reasonable for the considered case.

\begin{remark}  We point out that  the optimal choice of an agent at a time $s$ also depends on the congestion along the branches  (namely, the other components of $\rho$) that it will be running in the future. 
Nonetheless all statements would hold true also for a congestion function ${\cal F}^{w_1,w_2}$ depending on the whole vector $\rho$, instead that on $\rho^{w_1,w_2}$. Indeed, all along  the process of search of a fixed point, what is relevant is that the the vector $\rho$ is considered given, so that the values of the congestion costs are also given, and the arguments of the proofs contained in this paper would not change. 
\end{remark}

Note that, under the hypothesis (H3'),  the search of a fixed point described in Section 1 can be performed for $\rho$
rather than for $\mathcal{M}$. This is the subject of the next section.

%%%%%%%%%%%%%%%%%%%%%%%%%%%%%%%%%%%%%%%%%%%%%%%%%
%%%%%%%%%%%%%%%%%%%%%%%%%%%%%%%%%%%%%%%%%%%%%%%%%

\section{Existence of a mean field   equilibrium}

\label{sec:2} {In this section we give a proof of the existence of a mean
field equilibrium.}

Let $L(f)$ the Lipschitz constant of a function~$f$. As space to search for a fixed point, we choose
\begin{equation}\label{eq:X}
X=\left\{f:[0,T]\to\mathbb{R}_+:  {L(f)\le\tilde L},\  |f|\le K \right\}^4.
\end{equation}
the Cartesian product four times of the space of Lipschitzian functions with Lipschitz constant not greater than $\tilde L$  and overall bounded by $K$, where $\tilde L$ and $K$ are the constants  defined in Remark \ref{rholip}. 
Note that $X$ is convex and compact with respect to the uniform topology.

  We then search for a fixed point of the  {multi-function $%
\psi:X\to X$,  with
$\rho\mapsto \rho^\prime\in\psi(\rho)$ where the idea is to obtain 
$\rho^\prime$ as follows}: (i)  $\rho$ is inserted in  (\ref{eq:V00})--(\ref{eq:V11}),  the optimal control  {$u=-V_\theta$} is derived; (ii)  $u$ is inserted in (\ref{eq:transportXX}) and the distribution $m$ is derived; (iii)  $\rho^\prime$ is derived from (\ref{rho}). Actually, due to the hypothesis (H3'), in steps (ii) and (iii) we do not need to use the continuity equations for $m^{w_1,w_2}$ as in \cite{bagpes}:  instead from  the input   $\rho$  we may build the optimal controls and then the flow functions $g$ at the significant states, 
 as described in  Appendix  {A}, and finally the new mass concentrations $\rho'$ by means of (\ref{eq:totalmass}).

\begin{remark} Note that in general $\psi$ is \emph{not} single valued,  that is, $\psi(\rho)$ is a nonsigleton subset of $X$. Indeed the optimal control may not  be unique. For any fixed $t$, the minimization procedure in (\ref{eq:V00})--(\ref{eq:V11}) may return more than one minimizer. In particular,  {in} (\ref{eq:V01})--(\ref{eq:V11}) this  may happen even along   a whole time interval (on the contrary, in (\ref{eq:V00}) this may occur at most at a single instant). 
With  uniqueness of the optimal control  all agents in the same position at the same instant make the same choice (consistently with  mean field game models, where agents are homogeneous and indistinguishable).  When instead uniqueness fails,  $\psi$ can be built in many ways,  as many as the different optimal behaviors. \\
In \citet[Assumption 1]{bagpes}, the above difficulty was bypassed assuming  \textit{a-priori}   a finite number $N$  of times,  independent of $\rho$, at which those multiplicities appear. Here we drop such an assumption. \qed

\end{remark}

We will obtain $\psi $ and its fixed point $\bar\rho$ through a limiting procedure on a sequence $%
\left\{ \psi _{\varepsilon }\right\} _{\varepsilon >0}$ of functions
approximating $\psi $ in a suitable sense, and the corresponding fixed points $\rho_\varepsilon$. The single function 
$\psi _{\varepsilon }$ is obtained through~(i)--(iii) above, with
the difference that in (ii), rather than choosing optimal controls, one
chooses $\varepsilon -$optimal controls and, along time, an $\varepsilon-$\emph{optimal stream}.  We divide the construction into several steps.

\medskip

  \noindent \emph{Step 1: $\varepsilon$-optimal streams.}
\begin{definition} 
Assume the branch $B_{w_{1},w_{2}}$ is entered at 
 the state $(\widehat{\theta },w_{1},w_{2}),$ with $\widehat{\theta }\in \{\theta _{S},\theta _{1},\theta _{2}\}$,  and let $u^{w_{1},w_{2}}_{i}$, $i\in\{1,2,3\}$ be the controls defined through (\ref{oc1}) (\ref{oc2}) (\ref{oc3}) (\ref{oc4}). Consider also a partition  $\tau^{w_1,w_2}=\{t^n\}_n$ of the interval $[0,T]$, {and fix $\varepsilon>0$}. 
Then ${u}^{w_{1},w_{2}}_{\varepsilon }$ is an
  \emph{$\varepsilon -$optimal stream for $B_{w_1,w_2}$} associated to the partition $\tau^{w_1,w_2}$ if
$${u}^{w_{1},w_{2}}_{\varepsilon }(s)=u^{w_{1},w_{2}}_{i_n}(s), \ \ s\in[ t^{n},t^{n+1}[ $$  where  $u^{w_{1},w_{2}}_{i_n}$ is  optimal at $t^{n}$ and $\varepsilon$-optimal at all $s\in] t^{n},t^{n+1}[$, {that is, it realizes the minimum cost up to an error not greater than $\varepsilon$.}
 \end{definition}

An $\varepsilon -$optimal stream associated to a general partition $\tau$ may or may not exist, but it certainly does when the partition is refined enough.  Indeed, consider  the functions involved in  the minimization process in 
 (\ref{eq:V00}), (\ref{eq:V01}), (\ref{eq:V10}), (\ref{eq:V11}). 
If the minimum is realized up to the error $\varepsilon$, the minima   are attained within  intervals of type $[t+h_\varepsilon,T],$ where $h_\varepsilon$ is determined as $h$ in (\ref{h}), with due differences, so that the functions cited above are Lipschtiz continuous.
Denote by $L$ the greater of Lipschtiz constants of these functions. Then a control $u^{w_1,w_2}_i$ optimal at $t$ remain $\varepsilon$-optimal at least along $[t,t+\frac \varepsilon L[$.  

In addition note that, for every fixed $\varepsilon>0$, there may exist more than one $\varepsilon$-partition, as optimal controls at significant points may be multiple. Anyway, the number of $\varepsilon$-partitions is overall bounded by a number $M_\varepsilon$, as the optimal controls at every switching point are at most 3.
This argument proves the following Lemma.

\begin{lemma}\label{taueps}
Fixed $\varepsilon>0$, set $N=\max\{n\in{\mathbb N}: n<(TL)/\varepsilon\}+1$.  Consider the partition $\tau_\varepsilon$ of $[0,T]$ such that:
 \begin{itemize}
\item[$(i)$] $t^0_\varepsilon=0$; $t^N_\varepsilon=T$;
\item[$(ii)$] $t^n_\varepsilon=\frac{n L} \varepsilon$, for all $n\in\{1,\dots, {N-1}\}$.
\end{itemize}
Then the set of $\varepsilon$-optimal streams associated to $\tau_\varepsilon$ is nonempty and finite.
 \end{lemma}

 To better clarify the   definition of $\varepsilon$-optimal stream, we specify how the agents implementing it behave. For example, consider the state $(\theta_S,1,1)$ at time $t^0_\varepsilon=0$ and the $\varepsilon$-optimal stream on branch $B_{1,1}$, so that optimal controls that we refer to are chosen among those in (\ref{oc4}). 
Say,  the controls $u^{1,1}_{1}$ and  $u^{1,1}_{2}$  
are both optimal at $t_\varepsilon^0$. Then, at $t^{0}_\varepsilon$, agents all chose one of the two, say $u^{1,1}_{1}(t_\varepsilon^0)$. The agents arriving at $(\theta_S,1,1)$ at times $s > t^0_\varepsilon$ continue to choose the same  control $u^{1,1}_{1}(s)$  along the time interval $[0,L/\varepsilon]$ (where $u^{1,1}_{1}(s)$ is $\varepsilon$-optimal). 
Let $t_{\varepsilon }^{1}=\frac L\varepsilon$:
at this time, the agents arriving at $(\theta_S,1,1)$ update the strategy choice, selecting among  controls optimal  at $t^1_\varepsilon$ as before. Say  they choose $u^{1,1}_{2}(t_{\varepsilon }^{1})$. As before, agents arriving  at $(\theta_S,1,1)$ at time  $s$, with  $s > t^1_\varepsilon$,   continue to use   $u^{1,1}_{2}(s)$ until time $t^2_\varepsilon\equiv2\frac L\varepsilon$, knowing that  $u^{1,1}_{2}(s)$ is   $\varepsilon$-optimal on $[L/\varepsilon,(2L)/\varepsilon]$, and so on.

\begin{remark}  
{By definition, an agent  implementing a control chosen along a $\varepsilon-$optimal stream is using an $\varepsilon$-optimal control, that is, it is realizing the minimum of the payoff up to a maximum error of  $\varepsilon$. } 

Streams on different branches   all start at 0 and are defined  along the interval $[0,T]$. 
Note that an $\varepsilon$-optimal stream is not the strategy of a single agent at different times, whereas a strategy that is $\varepsilon$-optimal if adopted, along time, by the agents at $\hat\theta$, independently of the fact that any agent is present at $\hat\theta$ at that time. \color{black}
\qed
\end{remark}

\begin{remark} 
  Controls that are $\varepsilon $-optimal account for a
phenomenon known in literature as \emph{herd behavior} \citep
{Banerjee}: an agent is influenced by the decisions made by other agents,
even when that is not the optimal choice for him/her, unless the discrepancy
from the optimal choice is too large (greater than~$\varepsilon $).
\qed
\end{remark}

\medskip

  \noindent \emph{Step 2:  Split Fractions.}
For a fixed $\varepsilon>0$, consider the partition $\tau_\varepsilon$ defined in Lemma \ref{taueps}.  Consider then the vector
 $u_\varepsilon=(u_\varepsilon^{1,1},u_\varepsilon^{0,1},u_\varepsilon^{1,0},u_\varepsilon^{1;0,0},u_\varepsilon^{2;0,0})$ whose components are  $\varepsilon$-optimal streams associated to  $\tau_\varepsilon$ at the points $(\theta_S,1,1)$, $(\theta_1,0,1),$ $(\theta_2,1,0),$ $(\theta_1,0,0),$ $(\theta_2,0,0)$, respectively.

 As consequence of  Lemma \ref{taueps}, the $\varepsilon$-optimal vectors $u_\varepsilon$  are a finite number.

Note that once $u_\varepsilon $ is fixed,  all agents  choose equally, in time. Since we will need to consider the possibility for agents to split into fractions among different  vectors  $u_\varepsilon$  (that is, splitting among multiple optimal controls on  instants of $\tau_\varepsilon$), we give the following definition.

%At every point in $\tau_\varepsilon$ defined in Remark \ref{tau} the mass flow initiated by $g$  is split in fractions, according to the fraction of agents choosing each single available optimal controls at that point. Such fractions are described in the following definition.

\begin{definition}\label{def:split}Consider the partition $\tau_\varepsilon$ defined in Lemma \ref{taueps}. An \emph{$\varepsilon$-split function} is a vector  
 $\lambda_\varepsilon\in L^\infty(0,T)^{13}$, with coordinates
\begin{equation}
\label{eq:fractions}
\begin{array}{ll}
\displaystyle
\left(\lambda^{(\theta_S,1,1)}_1,\lambda^{(\theta_S,1,1)}_2,\lambda^{(\theta_S,1,1)}_3,\right.\\
\displaystyle
\lambda^{(\theta_1,0,1)}_1,\lambda^{(\theta_1,0,1)}_2,\lambda^{(\theta_1,0,1)}_3,
\lambda^{(\theta_2,1,0)}_1,\lambda^{(\theta_2,1,0)}_2,\lambda^{(\theta_2,1,0)}_3,\\
\displaystyle
\left.\lambda^{(\theta_1,0,0)}_1,\lambda^{(\theta_1,0,0)}_2,\lambda^{(\theta_2,0,0)}_1,\lambda^{(\theta_2,0,0)}_2\right).
\end{array}
\end{equation} 
constant on subintervals induced by the partition $\tau_\varepsilon$ and such that the \emph{split fractions} $\lambda^{(\hat \theta,w_1,w_2)}_i$ satisfy:
\begin{itemize}
\item[$(i)$]  $\lambda^{(\hat \theta,w_1,w_2)}_i(s)\ge 0,$ for all $s\in[0,T]$,
 and $\lambda^{(\hat \theta,w_1,w_2)}_i(t)=0$ if $u^{w_1,w_2}_i$ is not optimal at $(\hat \theta,w_1,w_2,t)$ for $t\in\tau_\varepsilon$;
\item[$(ii)$] $\sum_{i=1}^3\lambda^{(\hat \theta,w_1,w_2)}_i(s)=1$ for all $s\in[0,T]$.
\end{itemize}
\end{definition}

%%%%%%%%%%%%%%%%%%%%%%%%%%%%%%%%%%%%%%%%%%%
%%%%%%%%%%%%%%%%%%%%%%%%%%%%%%%%%%%%%%%%%%%
%%%%%%%%%%%%%%%%%%%%%%%%%%%%%%%%%%%%%%%%%%%
%%%%%%%%%%%%%%%%%%%%%%%%%%%%%%%%%%%%%%%%%%%
%%%%%%%%%%%%%%%%%%%%%%%%%%%%%%%%%%%%%%%%%%%
%%%%%%%%%%%%%%%%%%%%%%%%%%%%%%%%%%%%%%%%%%%

\bigskip

 \emph{Step 3: Construction of $\psi_\varepsilon(\rho)$.} 

Let $\varepsilon>0$ and $\rho=(\rho^{1,1},\rho^{0,1},\rho^{1,0},\rho^{0,0})\in X$ be fixed. Let also $\tau_\varepsilon$ be the partition of $[0,T]$ described above. 
We now build a multifunction
{$\psi_\varepsilon(\rho)\subseteq X$} with compact and convex images and closed
graph, to which later we can apply Kakutani fixed point theorem.

\begin{itemize}
\item[(a)] We define $\tilde\psi_\varepsilon(\rho)\subseteq X$ as the finite set of vectors $\rho^{\prime }=(\rho^{\prime 1,1},\rho^{\prime
		0,1};\rho^{\prime 1,0},\rho^{\prime 0,0})$ in $X$ constructed in the following
	way. For every significant point $(\hat\theta,w_1,w_2)$, we consider  an $\varepsilon-$optimal vector
	$u_{\varepsilon}$, associated to $\tau_\varepsilon$.
	Then, given the arrival flow of agents $g$, we use $u_{\varepsilon}$ to determine (forwardly from ($\theta_S,1,1)$, and at all significant points) all the flow functions $g_{w_1,w_2}$, as explained in   Appendix  {A}.
	Finally, we use $g_{w_1,w_2}$ and (\ref{eq:totalmass}) to compute the total mass~$\rho^{\prime }$. We repeat the construction  for all possible (and finite, by   Lemma \ref{taueps}) choices  $u_\varepsilon$ and call the set of all outcomes $\tilde\psi_\varepsilon(\rho)$. Note that $\tilde\psi_\varepsilon(\rho)$ is a finite set.

\item[(b)] We define $\psi_\varepsilon(\rho)\subseteq X$ as follows. We consider an arbitrary $\varepsilon$-split function $\lambda_\varepsilon$, and assume that at every point of $\tau_\varepsilon$ the incoming mass $g_{i,j}$ of agents split in fractions $\lambda^{(\hat\theta,w_1,w_2)}_1g_{i,j}$, $\lambda^{(\hat\theta,w_1,w_2)}_2g_{i,j}$, $\lambda^{(\hat\theta,w_1,w_2)}_3g_{i,j}$ (remaining at $\hat\theta$ or acceding the subsequent branches,  through the process explained in  Appendix  {A}). At the end of the process, the output is $\rho^\prime$.
We define $\psi_\varepsilon(\rho)$ as the set of all $\rho^\prime$  generated for all possible choices of the $\varepsilon$-split function $\lambda_\varepsilon$. Clearly $\psi_\varepsilon(\rho)\supset \tilde \psi_\varepsilon(\rho)$.

%   NON CANCELLARE
%\item[(b)]\added{We define  $\psi_\varepsilon(\rho)\subseteq X$ as the convex hull of $\tilde\psi_\varepsilon(\rho)$,
%	i.e., each quadruplet in $\rho'\in \psi_\varepsilon(\rho)$ is a convex combination of the quadruplets of $\tilde\psi_\varepsilon(\rho)$.
%	Specifically,  $\rho'\in \psi_\varepsilon(\rho)$ implies  $\rho' = \sum_{i\in I  } \alpha_i \rho'_i$, where $I$ is a finite set of indexes, $\rho'_i \in \tilde\psi_\varepsilon(\rho)$ and $\alpha_i $ are nonnegative coefficient such that  $\sum_{i\in I} \alpha_i = 1$. 
%	In addition, let $u_{\varepsilon,i}$ the $\varepsilon-$suboptimal control that generates $ \rho'_i$, for all $  \rho'_i \in \tilde\psi_\varepsilon(\rho)$.
%	Due  to the linearity of the operations described in the previous point (a), 
%	$\rho'$ can be seen as generated by splitting the arrival flow of agents $g$ in fractions proportional to $\alpha_i$, 
%	each of which uses a control $ u_{\varepsilon,i}$, which in turn generates the total mass~$\alpha_i \rho'_i$.}
%
%    NON CANCELLARE 

\end{itemize}
%
%$g_{w_1,w_1}$
% and the functions  described above and its convex hull $\Psi$. The following Lemma holds true.

\begin{lemma}\label{lem:2}  The set $\psi_\varepsilon(\rho)$ is a {non-empty} convex and compact subset of $X$, {for any $\rho \in X$}.  
In addition, the map $\rho\mapsto \psi_\varepsilon(\rho)$ has closed graph and, in particular, 
it has a fixed point $\rho_\varepsilon\in X$.
\end{lemma} 

\begin{proof} We preliminarily observe that the set $\psi_\varepsilon(\rho)$ is  non-empty by construction. In addition, it is also closed, and hence, since $X$ is compact, it is compact. The closedness comes from the fact that the multifunction $\rho\mapsto\psi_\varepsilon(\rho)$ has closed graph, which will be proved below.
To prove that $\psi_\varepsilon(\rho)$ is convex, we select $\xi,\eta\in\psi_\varepsilon(\rho)$, and $\alpha\in[0,1]$ and show that also $\alpha \xi+(1-\alpha)\eta
\in\psi_\varepsilon(\rho)$. Assume that $\lambda_\varepsilon^\xi$ and $\lambda_\varepsilon^\eta$ are the splits generating $\xi$ and $\eta$ respectively. Then $\alpha \xi+(1-\alpha)\eta$ is generated by $\alpha\lambda_\varepsilon^\xi+(1-\alpha)\lambda_\varepsilon^\eta$, which is itself an $\varepsilon$-split function, and the proof is complete.

Next, prove that the map has closed graph.
We consider a sequence $\{\rho^n\}\subset X$, with  $\rho^n\to\rho$ in $X$ 
	and we prove that for every selection $\rho^{\prime n}\in\psi_\varepsilon(\rho^n)$, 
	with $\rho^{\prime n}\to\rho^{\prime}$ in $X$, we have $\rho^{\prime}\in\psi_\varepsilon(\rho)$. We divide the long proof into several steps. Note that all along the partition $\tau_\varepsilon=\{t^n_\varepsilon\}$ is the same for all $\rho^n, \rho$.

(1) Consider the the value functions $V^n$ and 
$V$, defined  by (\ref{eq:V00}),(\ref{eq:V01}),(\ref{eq:V10}),(\ref{eq:V11}) and associated, respectively, to the  choices of masses $\rho^n$ and $\rho$  in the congestion cost
$\cal F$. By Theorem \ref{th:valuefun} and hypothesis (H3') they are equi-bounded and equi-Lipschitz in time, and continuously depending on $\rho^n$, $\rho$ respectively. Since $\rho^n\to\rho$ uniformly in time, then 
 $V^n\to V$ uniformly on $[0,T]$, and   the corresponding minimizers $\tau^n, \eta^n\in[0,T]$ of the arguments in (\ref{eq:V00}),(\ref{eq:V01}),(\ref{eq:V10}),(\ref{eq:V11}) are converging (at least along a subsequence) to suitable $\tau$ and $\eta$, minimizing the same arguments with  $\rho$ in place of $\rho^n$. 

Then, along the same subsequence, at any instant $t_\varepsilon^n$ of the partition $\tau_\varepsilon$,  the optimal controls $u^n_\varepsilon(t)$ induced by $\rho^n$ converge to an optimal control $u_\varepsilon(t)$  induced by $\rho$.

(2)  Every $\rho^{\prime n}$ is generated by some $\varepsilon$-split functions $\lambda_\varepsilon^n$ associated to $\rho^n$. We prove next  that $\rho'$ is also generated by a split function $\lambda_\varepsilon$ associated to $\rho$. 

Since     $\varepsilon$-optimal vectors $u_\varepsilon$ are finite, we may assume that  (possibly along a subsequence) in every subinterval $[t_\varepsilon^n,t_\varepsilon^{n+1}[$ the active components of all   $\varepsilon$-split functions $\lambda_\varepsilon$ are the same:
if a component of $\lambda_\varepsilon^m$ is nonzero in a subinterval, then the same component  is nonzero for all other $\varepsilon$-split functions $\lambda_\varepsilon^n$.  Moreover, since   $\varepsilon$-split functions are constant in all subintervals, we can also suppose that the sequence $\lambda_\varepsilon^n$ (uniformly) converges to a limit function $\lambda_\varepsilon$ which is also a $\varepsilon$-split function.

(3) We examine first the branch $B_{1,1}$ and assume w.l.o.g.   that along $[0,t^1_\varepsilon]$ the total mass of entering agents is nonzero (i.e. the integral of the given flow $g$ is not null). Note that $\rho^{\prime n}$ is uniformly converging to $\rho'$, and that branches are initially empty (see hypothesis (H2)). These facts would be in contradiction with a first component of $\rho'$ (i.e. the component on $B_{1,1}$)   not generated by the limiting process described in (1) and (2).
Indeed, the second and the third components of $\rho^{\prime n}$ and $\rho'$ on the branches $B_{0,1}$ and $B_{1,0}$ respectively,  are given by the flow functions as in (\ref{eq:totalmass}), which strongly depend on the optimal controls and on the split functions on $B_{1,1}$ (see  Appendix A).

(4) We iterate the argument in (3) both for subsequent intervals of $\tau_\varepsilon$ and for the other branches,  finally obtaining  the closed graph property of $\psi_\varepsilon$. 

Hence, by Kakutani fixed point theorem,
 the map $\rho\mapsto \psi_\varepsilon(\rho)$  has a fixed point.\qed

\end{proof}

{Hereinafter, we denote by $\rho_\varepsilon$ a fixed point for $\psi_{\varepsilon}(\rho)$, i.e., a total mass that satisfies $\rho_\varepsilon \in \psi_{\varepsilon}(\rho_\varepsilon)$}.

%\begin{remark} An interpretation of the above lemma and the definition of $\rho_\varepsilon$ is the following.
%Assume that all agents implement a $\varepsilon-$optimal stream and predict that the realized congestion will be given by $\rho_\varepsilon$. 
%There is a  choice among the available  $\varepsilon-$streams  such
%that $\rho_\varepsilon$ is actually realized and, moreover,
%every agent is subject to a cost (\ref{eq:cost}) whose difference with the
%optimal one not greater than $\varepsilon$.\qed
%\end{remark}

\medskip

Before stating the existence of a mean field   equilibrium, 
we introduce the following definitions that help restrict the equilibrium concept to the purpose of our problem.

\begin{definition}
\label{def:sf}  
A {\it split function} 
 is a vector   $\lambda\in L^\infty(0,T)^{13}$, with components given by (\ref{eq:fractions}) such that $\lambda^{(\hat \theta,w_1,w_2)}_i$ satisfy:
\begin{itemize}
\item[$(i)$]  $\lambda^{(\hat \theta,w_1,w_2)}_i(s)\ge 0,$ for a.a. $s\in[0,T]$,
 and $\lambda^{(\hat \theta,w_1,w_2)}_i(s)=0$ if $u^{w_1,w_2}_i$ is not optimal at $s$;
\item[$(ii)$] $\sum_{i=1}^3\lambda^{(\hat \theta,w_1,w_2)}_i(s)=1$ for a.a. $s\in[0,T]$.
\end{itemize}

\end{definition}

Note that the definition differs from Definition \ref{def:split} in the fact that split function are not linked to any $\varepsilon$-partition, and they are not necessarily piecewise constant.
\begin{definition}
\label{def:MFGE}
Let $\psi$ and $\psi_\varepsilon$ be the functions described at the beginning of Section \ref{sec:2}.

A \textit{$\varepsilon$-mean field   equilibrium} is a total mass~$\rho_\varepsilon \in X$ that satisfies $\rho_\varepsilon \in \psi_{\varepsilon}(\rho_\varepsilon)$. 

A \textit{mean field   equilibrium} is a total mass $\rho\in X$ that satisfies $\rho \in \psi(\rho)$.

\end{definition}

 Note that   $\rho \in \psi(\rho)$ implies 
that $\rho$ induces a set of optimal controls as in~(\ref{oc1}),~(\ref{oc2}),~(\ref{oc3}),~(\ref{oc4}), used by masses of agents fractioned according to a split function $\lambda$ in every branch: the flows $\lambda^{(\theta_S,1,1)}_i\;g$ generate the corresponding flows entering in the branches $B_{0,1},B_{1,0}$;  incoming flows in branches $B_{0,1},B_{1,0}$ 
 split again according to   $\lambda$, and generate flows entering   $B_{0,0}$, so that the final outcome, ruled  by (\ref{eq:totalmass}), is again $\rho$.

\begin{theorem}\label{th:MFE}
	Assume  \emph{(H1),(H3'),(H4)}. Then there exists a mean field   equilibrium.
\end{theorem}

\textit{Proof.} Consider $\varepsilon>0$, the $\varepsilon$-partition $\tau_\varepsilon$, a fixed point $%
\rho_\varepsilon$ for $\psi_\varepsilon$, and  the associated  $\varepsilon$-streams $ u_\varepsilon $ and  
 $\varepsilon$-split function $\lambda_\varepsilon\in
L^\infty(0,T)^{13}$. Recall that $\lambda_\varepsilon$
is constant along subintervals induced by $\tau_\varepsilon$. (Note also  that the number of subintervals is not uniformly bounded with respect to $\varepsilon$). Arguing as in  (3) in Lemma \ref{lem:2} we see that, at least along a subsequence, $\lambda_\varepsilon$
converges in the weak star topology to a vector $\lambda\in L^\infty(0,T)^{13}$, as $\varepsilon\to0$, that is
$\int_0^T\lambda_\varepsilon \mu\,dt\overset{\varepsilon\to0}{\to}\int_0^T\lambda \mu\,dt$, for all 
function $\mu\in L^1(0,T)$, and $\rho_ \varepsilon$ converges uniformly to $\rho\in X$. 
Hence by  construction   $\lambda$ is a split  function generating the total mass $\rho$.
%\addedf{It generates flow functions $g_{\cdot,\cdot}$ which by (\ref{eq:totalmass})
%generate $\rho$. 
\footnote{Note that (\ref{eq:totalmass}) are integral equalities, and  that every flow function $g_{i,j}$ is built by means of the incoming flow $g$, the optimal controls and the split functions, as explained in  Appendix A. Hence the use weak star convergence of the split functions is appropriate.} 
What is left to show is that whenever a component of $\lambda$ is not null, then the corresponding choice of the control is optimal.
 For almost every instant $t$ such that one of the components of $%
\lambda$ is not null, by   weak-star convergence there exists
a sequence  of instants $t_\varepsilon$  converging to $t$
such
that, at least along a subsequence, the corresponding component of $\lambda_\varepsilon(t_\varepsilon)$ is
not null.

Indeed, if this is not true, then there exists a neighborhood of $t$ where, for $\varepsilon$ sufficiently small, 
that component of all $\lambda_\varepsilon$ are null a.e.. By the weak star convergence, that would mean  that the corresponding component of $\lambda$ is also null a.e. on that neighborhood, which is a contradiction.

This means that, for every $\varepsilon$ the  choice of the corresponding component of $u_\varepsilon$   is $\varepsilon$-optimal at $t_\varepsilon$. 
 The conclusion is then reached  sending $\varepsilon\to0$ and arguing as in (1) in the proof of Lemma \ref{lem:2}. \qed

%\begin{table}
%\begin{tabular}{lll}
%\hline\noalign{\smallskip}
%\end{tabular}
%\end{table}

%\begin{acknowledgements}
%\end{acknowledgements}



\section{The optimization problem}

In this section we describe an optimization problem faced by a
local authority, that we refer to as \emph{controller}, and that needs to manage the tourists' flow in a city.
As mentioned in the Introduction, this problem can be framed in the new class of mean field games proposed by  \cite{LLions2018}.

We restrict our analysis to congestion cost functions of the form
\begin{equation}  \label{eq:F}
\mathcal{F}^{w_1,w_2}(\rho)=\alpha_{w_1,w_2}\rho^{w_1,w_2}(s)+\beta_{w_1,w_2}
\end{equation}
{with $\rho=(\rho^{1,1},\rho^{0,1},\rho^{1,0},\rho^{0,0})\in X$.} In~(\ref%
{eq:F}), the coefficients $(\alpha_{w_1,w_2},\beta_{w_1,w_2})$ are chosen by the controller, aiming 
 to force the equilibrium  to be as close as possible (in uniform
topology) to a reference string $\overline\rho\in X$, i.e. to minimize: 
\begin{equation}  \label{eq:distance0}
\max_{w_1,w_2\in\{0,1\}}\left\{\max_{t\in[0,T]}|\rho^{w_1,w_2}(t)-\overline%
\rho^{w_1,w_2}(t)|\right\}=\|\rho-\overline\rho\|_X.
\end{equation}
Let us denote  by $\chi_{\alpha,\beta}$ the set of 
mean field game equilibria corresponding to the choice of
  parameters $\alpha=(\alpha_{1,1},\alpha_{0,1},\alpha_{1,0},\alpha_{0,0})$
and $\beta=(\beta_{1,1},\beta_{0,1},\beta_{1,0},\beta_{0,0})$, which are assumed to belong to a compact set $\mathcal{K}\subset\mathbb{R}^4\times\mathbb{R}^4$. 
Then the controller faces the
optimization problem given by
\begin{equation}  \label{eq:distance1}
\inf_{(\alpha,\beta)\in \mathcal{K}}\mathcal{E}(\alpha,\beta),
\end{equation}
where $\mathcal{E}(\alpha,\beta) = \inf_{\rho\in\chi_{\alpha,\beta}}\|\rho-\overline\rho\|_X$.

\begin{theorem}
	In the same assumptions of Theorem~\ref{th:MFE}, there exists an optimal pair $(\alpha,\beta)\in \mathcal{K}$ that solves problem~(\ref{eq:distance1}).
\end{theorem}

\textit{Proof.} Let $(\alpha^n,\beta^n)\in \mathcal{K}$ be a minimizing sequence for $%
\mathcal{E}$, and for every $n$ let $\rho^n\in\chi_{\alpha_n,\beta_n}$
realize the infimum in (\ref{eq:distance1}) up to an error  $1/n$. By compactness of $\mathcal K$ and $X$, we
may suppose that $(\alpha^n,\beta^n)$ converges to $(\tilde\alpha,\tilde%
\beta)\in \mathcal{K}$, and that $\rho^n$ uniformly converges to $\tilde\rho\in X$. Hence we
 only need to prove that $\tilde\rho\in\chi_{\tilde\alpha,\tilde\beta}$, 
that is, $\tilde\rho$ is a mean field   equilibrium. 
Let $\lambda^n$ be the split function (Definition \ref{def:sf})  associated to the mean field 
equilibrium $\rho^n$, and let $\tilde\lambda$ be its weak-star limit. Arguing as in the proof of  Theorem~\ref{th:MFE}, we
derive that $\lambda$ is a split  function generating the total mass~$\tilde\rho$, and hence $\tilde\rho$ is a mean field equilibrium
for $(\tilde\alpha,\tilde\beta)$.
\qed

\begin{remark}
	A   variation of problem~(\ref{eq:distance1}) is the following
	
	\[
	\mathcal{E}(\alpha,\beta) = \sup_{\rho\in\chi_{\alpha,\beta}}\|\rho-\overline\rho\|_X
	\]
	
	\noindent
	The interest of this new problem relies on the fact that  the controller tries to manage the  worst case scenario. 
	Unfortunately, in this case the existence of an optimal pair is not evident. 
	It can be proved that there exists a pair $(\alpha,\beta)$ and $\rho\in{\cal X}_{\alpha,\beta}$ such that $\|\rho-\overline\rho\|_X$ is equal to the infimum of $\mathcal{E}(\alpha,\beta)$. However, we are not able to prove that  $(\alpha,\beta)$ is optimal, due to the possible presence of multiple mean field equilibria. 

	A  way to  bypass  the above difficulty could be to  guarantee the uniqueness of the equilibrium,
	for example assuming stronger hypotheses on the cost $\cal F$. We leave such matter for future investigations. \qed
\end{remark}

\medskip

	So far we have assumed that the coefficients $(\alpha_{w_1,w_2},%
	\beta_{w_1,w_2})$ of congestion cost functions $\mathcal{F}^{w_1,w_2}$ are
	constant over time.  
	We now generalize  by assuming that $(\alpha_{w_1,w_2},\beta_{w_1,w_2})$   piecewise constant   over time and
that the control is implemented at the significant points of each branch of
the network, e.g. through gates. 

We consider a finite sequence of fixed instants $t_0=0<t_1<t_2<...<t_N=T$
and, for every $i=0,..N-1${, the coefficients $(\alpha^i,\beta^i)\in \mathcal{K}$ chosen by the
controller in  $[t_i,t_{i+1}[$.} {The
congestion cost paid by an agent at time $s$ becomes 
\begin{equation*}
{\mathcal{F}}^{w_1,w_2}(s)=\alpha_{w_1,w_2}(\tau(s))\rho^{w_1,w_2}(s)+%
\beta_{w_1,w_2}(\tau(s))
\end{equation*}
where } {$\alpha_{w_1,w_2}(\tau)=\alpha_{w_1,w_2}^i$, respectively $%
\beta_{w_1,w_2}(\tau)=\beta_{w_1,w_2}^i$, for $t_i\le \tau<t_{i+1}$, $%
i=0,...,N-1$; and $\tau(s)\le s$ is  % the last switching instant not greater than $s $ along the agent trajectory, i.e., 
the instant at which the agent state
entered $B_{w_1,w_2}$.}

{In this situation the total} cost payed by an agent is  the usual 
\begin{equation}  \label{eq:cost2}
\begin{array}{ll}
J(u;\theta,t,w_1,w_2)=\int_t^T\left(\frac{u(s)^2}{2}+{\mathcal{F}}^{w_1,w_2}(s)%
\right)ds+ &  \\ 
+c_1w_1(T)+c_2w_2(T)+c_S\xi_{\theta=\theta_S}(T) & 
\end{array}%
\end{equation}
but  the congestion cost explicitly depends on   time 
through $\alpha_{w_1,w_2}$ and $\beta_{w_1,w_2}$.  Nevertheless,
such a dependence is compatible with the structure of the choices of the
agents in our model as presented in the following example. 

The agents that
arrive at $(\theta_S,1,1)$ at time $t$  they are given the
parameters $\alpha_{1,1}(t)$ and $\beta_{1,1}(t)$, which remain constant
until their exit from the {branch} $B_{1,1}$.
 Suppose that the agents now
move to the branch~$B_{0,1}$. Then, at that moment, they are given the new values
of  parameters~$\alpha_{0,1}$ and $\beta_{0,1}$ which, again remain
 constant until their exit from the branch, and so on. 
This behaviour is justified as we assume that the controller applies its
controls at the beginning of the network branches and that 
the agents make their decisions when they enter a new branch: then a possible change in the parameters $%
(\alpha^i,\beta^i)$ is not perceived by the agents  ``on the road"
along a branch. {Consequently, the considered} time-dependence of the cost
does not change the arguments in previous sections. In particular, formulas (\ref{eq:V00})--(\ref{eq:V11}) remain
true, and the existence of a mean field   equilibrium is guaranteed.
Finally, the optimization problem (\ref{eq:distance1}) can be successfully
solved.

\section{Conclusions}

In this paper we have introduced a mean field game model representing   the flows of tourists in the alleys of a heritage city.
We then used this model to solve an optimization problem where the controller is a local authority
aiming to manage  tourists' flow in the direction of a target congestion of the city.

Further refinements of both the model and the optimization problem
  are possible from both a theoretical and an applicative
perspective. As an example, the assumption of the existence of surveillance
cameras, to count people entering and leaving areas of interest, may suggest
variations of the model and justify the
definition of a feedback control policy on the control parameters $%
(\alpha,\beta)$.
In addition one may consider   different objectives of the controller, or  the generalization 
of the model to a number $n$ of sites
 of interest. Note that, although 
in the latter case the number of alternative paths would increase
exponentially,  the direct experience of the authors in Venice
suggests that the vast majority of   tourists is interested in very few
attractions in a city.

\begin{acknowledgements}
This research was partially funded by a INdAM-GNAMPA project 2017
\end{acknowledgements}

\bibliographystyle{spbasic}      % basic style, author-year citations
\bibliography{Bagagiolo}

\begin{thebibliography}{28}
\providecommand{\natexlab}[1]{#1}
\providecommand{\url}[1]{{#1}}
\providecommand{\urlprefix}{URL }
\expandafter\ifx\csname urlstyle\endcsname\relax
  \providecommand{\doi}[1]{DOI~\discretionary{}{}{}#1}\else
  \providecommand{\doi}{DOI~\discretionary{}{}{}\begingroup
  \urlstyle{rm}\Url}\fi
\providecommand{\eprint}[2][]{\url{#2}}

\bibitem[{Achdou et~al.(2012)Achdou, Camilli, and Capuzzo~Dolcetta}]{AcCamDolc}
Achdou Y, Camilli F, Capuzzo~Dolcetta I (2012) {Mean field games: numerical
  methods for the planning problem}. {SIAM J Control Optim} 50:77--109

\bibitem[{Ambrosio et~al.(2008)Ambrosio, Gigli, and Savar\`e}]{ambgigsav}
Ambrosio L, Gigli N, Savar\`e G (2008) {Gradient Flows in Metric Spaces and in
  the Space of Probability Measures}. Lectures in Mathematics, Birkh\"auser
  Verlag, Basel, CH

\bibitem[{Bagagiolo and Danieli(2012)}]{bagdan}
Bagagiolo F, Danieli K (2012) Infinite horizon optimal control problems with
  multiple thermostatic hybrid dynamics. {Nonlinear Anal Hybrid Syst} 6:824 --
  838

\bibitem[{Bagagiolo and Pesenti(2017)}]{bagpes}
Bagagiolo F, Pesenti R (2017) Non-memoryless pedestrian flow in a crowded
  environment with target sets. In: Apaloo J, Viscolani B (eds) {Advances in
  Dynamic and Mean Field Games. ISDG 2016. Annals of the International Society
  of Dynamic Games}, 15, Springer, Cham, CH, pp 3--25

\bibitem[{Bagagiolo et~al.(2017)Bagagiolo, Bauso, Maggistro, and
  Zoppello}]{BBMZop}
Bagagiolo F, Bauso D, Maggistro R, Zoppello M (2017) {Game theoretic
  decentralized feedback controls in Markov jump processes}. {J Optim Theory
  Appl} 173:704--726

\bibitem[{Banerjee(1992)}]{Banerjee}
Banerjee AV (1992) {A Simple Model of Herd Behavior}. {Q J Econ} 107:797--817

\bibitem[{``Brian"~Park et~al.(2009)``Brian"~Park, Yun, and Ahn}]{BrianYun}
``Brian"~Park B, Yun I, Ahn K (2009) Stochastic optimization for sustainable
  traffic signal control. {Int J Sustain Transp} 3:263--284

\bibitem[{Camilli et~al.(2015)Camilli, Carlini, and Marchi}]{CCMar}
Camilli F, Carlini E, Marchi C (2015) A model problem for mean field games on
  networks. Discrete Contin Dyn Syst 35:4173--4192

\bibitem[{Camilli et~al.(2017)Camilli, De~Maio, and Tosin}]{CDeMTos}
Camilli F, De~Maio R, Tosin A (2017) Transport of measures on networks. Netw
  Heterog Media 12:191--215

\bibitem[{Carlier and Santambrogio(2012)}]{CarlierSant}
Carlier G, Santambrogio F (2012) A continuous theory of traffic congestion and
  wardrop equilibria. J Math Sci 181:792--804

\bibitem[{Carlier et~al.(2008)Carlier, Jimenez, and Santambrogio}]{Carlier1}
Carlier G, Jimenez C, Santambrogio F (2008) Optimal transportation with traffic
  congestion and wardrop equilibria. SIAM J Control Optim and 47 pp 1330--1350

\bibitem[{Cole et~al.(2006)Cole, Dodis, and Roughgarden}]{Cole}
Cole R, Dodis Y, Roughgarden T (2006) How much can taxes help selfish routing?
  J Comput Syst Sci 72:444--467

\bibitem[{Como et~al.(2013)Como, Savla, Acemoglu, Dahleh, and
  Frazzoli}]{ComoSavla}
Como G, Savla K, Acemoglu D, Dahleh MA, Frazzoli E (2013) Stability analysis of
  transportation networks with multiscale driver decisions. SIAM J Control
  Optim 51:230--252

\bibitem[{Cristiani et~al.(2015)Cristiani, Priuli, and Tosin}]{CPTos}
Cristiani E, Priuli F, Tosin A (2015) Modeling rationality to control
  self-organization of crowds. SIAM J Appl Math 75:605--629

\bibitem[{Dial(1999)}]{Dial}
Dial RB (1999) {Network-optimized road pricing: Part I: A parable and a model}.
  Oper Res 47:54--64

\bibitem[{Gomes and Horowitz(2006)}]{GomesHoro}
Gomes G, Horowitz R (2006) Optimal freeway ramp metering using the asymmetric
  cell transmission model. Transport Res C-Emer 14:244--262

\bibitem[{Hegyi et~al.(2005)Hegyi, De~Schutter, and
  Hellendoorn}]{HegyiSchutter1}
Hegyi A, De~Schutter B, Hellendoorn H (2005) Model predictive control for
  optimal coordination of ramp metering and variable speed limits. Transport
  Res C-Emer 13:185--209

\bibitem[{Hoogendoorn and Bovy(2005)}]{Hoogendoorn}
Hoogendoorn SP, Bovy PL (2005) Pedestrian travel behavior modeling. Netw Spat
  Econ 5:193–216

\bibitem[{Huang et~al.(2006)Huang, Caines, and Malhame}]{HCAINMAL}
Huang M, Caines P, Malhame R (2006) {Large population dynamic games:closed-loop
  McKean- Vlasov systems and the Nash certainly equivalence principle}. Commun
  Inf Sys 6:221--252

\bibitem[{Khattak et~al.(1996)Khattak, Polydoropoulou, and Ben-Akiva}]{Khattak}
Khattak A, Polydoropoulou A, Ben-Akiva M (1996) Modeling revealed and stated
  pretrip travel response to advanced traveler information systems. Transport
  Res Rec 1537:46--54

\bibitem[{Lachapelle et~al.(2010)Lachapelle, Salomon, and
  Turinici}]{LachSalTur}
Lachapelle A, Salomon J, Turinici G (2010) Computation of mean-field equilibria
  in economics. Math Models Methods Appl Sci 20:1--22

\bibitem[{Lasry and Lions(2006)}]{LLions}
Lasry J, Lions P (2006) Jeux \`{a}champ moyen ii {H}orizon fini et controle
  optimal. C R Math 343:679--684

\bibitem[{Lasry and Lions(2018)}]{LLions2018}
Lasry J, Lions P (2018) Mean-field games with a major player. C R Math
  356:886--890

\bibitem[{Morrison(1986)}]{Morrison}
Morrison S (1986) A survey of road pricing. Transport Res A-Pol 20:87--97

\bibitem[{Pigou(1920)}]{Pigou}
Pigou AC (1920) {The Economics of Welfare New York}. Macmillan, Trenton, NJ

\bibitem[{Smith(1979)}]{Smith}
Smith MJ (1979) The marginal cost taxation of a transportation network.
  Transport Res B-Meth 13:237--242

\bibitem[{Srinivasan and Mahmassani(2000)}]{Mahamassani}
Srinivasan K, Mahmassani H (2000) Modeling inertia and compliance mechanisms in
  route choice behavior under real-time information. Transport Res Rec
  1725:45--53

\bibitem[{Wardrop(1952)}]{Wardrop}
Wardrop JG (1952) Some theoretical aspects of road traffic research. Proc Inst
  Civ Eng 1:325--362

\end{thebibliography}

\appendix 

\section{Appendix: on functions $g$ and $g_{i,j}$}

{In this appendix, we discuss the relationship between the flow functions $g$
and $g_{i,j}$ and the components of  $\mathcal{M}$, and hence of $\rho$ both from a spatial and a temporal perspective.}
We consider only the branch $B_{1,1}$ as analogous arguments apply to the other branches of our model.

\smallskip

The agents'  flow, represented by the time depending $g$, enters the branch $B_{1,1}$ at $\theta_S$. Then $g$ can be interpreted as number of agents per unit of time and represents an incoming \emph{density with respect to time}.
Differently,  the total mass~$\rho^{1,1}$ on branch $B_{1,1}$ using~(\ref{rho}), 
is obtained  integrating  over the branch the \emph{density with respect to space} $m^{1,1}$, representing  the number of agents per space unit. The mass spread along the branch reaches in time
the  switching points $\theta_1,\theta_2$, where it is again converted in time dependent flows  $g_{i,j}$ entering the subsequent branches, with  $g_{i,j}$ again densities with respect to time.  The process repeats at every switching point.

We then describe  the mathematical relationship between these 
two different kinds of densities.
Our argument is connected  to what is called \emph{disintegration of a measure} 
\citep{ambgigsav,CDeMTos}.

Let us consider an agent at the position $\theta$ at time $s$ and assume that it is going towards~$\theta_1$. 
It arrived in $\theta_S$ at time $\underline \sigma(\theta,s) \leq s$  and it will reach $\theta_1$ at a time $\overline \tau(\theta,s) \geq s$. 
By Remark \ref{ccp}, in $s$ and in all of the interval $[\underline \sigma, \overline \tau)$ the agent has a constant velocity $u^{1,1}(s) =  u^{1,1}(\underline \sigma)$ such that $|u^{1,1}(\underline \sigma)|= |\theta_1-\theta_S|/(\overline \tau-\underline \sigma)$, 
until leaving the branch.
When the total mass $\rho$ is given, %(i.e. when performing the fixedpoint procedure), 
 then $\overline \tau$ realizes the minimum in
 $\inf_{\tau\in(\underline \sigma,T]}\left\{\frac{1}{2}\frac{(\theta _{1}-\theta _{S})^{2}}{\tau -\underline \sigma}
+\int_{t}^{\tau }\mathcal{F}^{1,1}\,ds +V(\theta _{1},\tau ,0,1)\right\}$, i.e., the second argument of the minimum in~(\ref{eq:V11}); 
the value of $\underline \sigma$ is evaluable using conditions
 (\ref{eq:V00})--(\ref{eq:V11}). 
 In addition, $\underline\sigma$ enjoys the properties linking $g$ and $g_{i,j}$ which are described below.

 Let us denote $q(\underline \sigma)=1/u^{1,1}(\underline\sigma)$ and,
 for the sake of simplicity, let us suppose that $\underline\sigma$ and $q$ are
differentiable (actually,  they are Lipschitz and
  at least  a.e. differentiable).

Let us also assume  $\theta_S=0$, $\theta>0$. We can implicitly define the value of  $\underline\sigma(\theta,s)$ through the following equation: 
\begin{equation}\label{eq:cambiodival}
s = \underline\sigma(\theta,s)+q(\underline\sigma(\theta,s))\theta.
\end{equation}
Differentiating (\ref{eq:cambiodival}) with respect to $s$ and $\theta$ we derive
\begin{equation*}
\underline\sigma_s(\theta,s)=\frac{1}{1+q^{\prime }(\underline\sigma(\theta,s))\theta},  \quad \underline\sigma_\theta(\theta,s)=\frac{q(\underline\sigma(\theta,s))}{1+q^{\prime }(\underline\sigma(\theta,s))\theta} \quad \Rightarrow \quad
\underline\sigma_s(\theta,s)=-\frac{\underline\sigma_\theta(\theta,s)}{q(\underline\sigma(\theta,s))}
\end{equation*}
Now, {suppose that all agents entering $(\theta_s,1,1)$ at any time move towards $\theta_1$. Then the flow $g$ of agents arriving in $\theta_S$ and moving
towards $\theta_1$ is spread over  $B_{1,1}$,
 according to the law $m^{1,1}(\theta,s)=-g(\underline\sigma(\theta,s))\underline\sigma_\theta(\theta,s)$.
In addition, the flow of agents crossing   $\theta$ in $B_{1,1}$ at time $s$ is given by $m^{1,1}(\theta,s)u^{1,1}(s) = g(\underline\sigma(\theta,s))\underline\sigma_s(\theta,s)$.
 Both relationships may be verified by standard mass balance/conservation arguments.
They obviously hold only if at time~$s$   agents have already arrived at $\theta$, that is when $\underline\sigma(\theta,s) \geq 0$, otherwise the density $m^{1,1}(\theta,s)$ is equal to zero. 
In particular, at the switching point $\theta_1$ the arriving flow, coinciding with the  flow $g_{0,1}(s)$ entering the new  branch $B_{0,1}$ in time,  is $s\mapsto
g_{0,1}(s) = g(\underline\sigma(\theta_1,s))\underline\sigma_s(\theta_1,s)$.

If differently, agents entering through $(\theta_s,1,1)$ split among different choices, and the
corresponding split fraction that moves towards $\theta_1$ is $\lambda_2^{(\theta_S,1,1)}$ (see (\ref{eq:fractions}) and Definition \ref{def:sf}), then the entering flow in $B_{0,1}$ through $\theta_1$ is $
s\mapsto g_{0,1}(\underline\sigma(\theta_1,s))=\lambda_2^{(\theta_S,1,1)}g(\underline\sigma(\theta_1,s))\underline\sigma_s(\theta_1,s)$. This is also the flow  $g_{01}$ to be considered in (\ref{eq:totalmass}).
Similar considerations (with different function $\underline\sigma$ and $q$)
hold in the case of agents moving towards $\theta_2$ in the branch $B_{1,1}$
and for all other cases in the other branches (with the corresponding flows $g_{i,j}$).

\bigskip
Let us finally argue on the uniqueness of   $\underline \sigma$ and $\overline \tau$. Specifically, consider the agents entering at time $\underline\sigma$ at $(\theta_S,1,1)$ and moving towards $\theta_1$, and that  reach such state  at time  $\overline \tau$. We claim: (1) any arrival time $\overline \tau$ originates from a unique entering time $\underline\sigma$;  
(2) any entering time $\underline \sigma$ generates a unique arrival time $\overline\tau$, for any nonzero flow of agents at $\theta_1$.  A sketch of the proofs follows:
 
(1) Let us suppose that the agents are optimally moving from $\theta_S$ to $\theta_1$ in the branch $B_{1,1,}$, and that the agents respectively starting at $\underline\sigma_1$ and at $\underline\sigma_2>\underline\sigma_1$ reach $\theta_1$ at the same time $\overline\tau<T$. This means that $\overline\tau$ optimizes the second term inside the minimum in (\ref{eq:V11}) for both $t=\underline\sigma_1$ and   $t=\underline\sigma_2$. Suppose that the function $\tau\mapsto V(\theta_S,\tau,1,1)$ is differentiable at $\overline\tau$. Since $\overline\tau$ is interior to $]\underline\sigma_2,T[\subset]\underline\sigma_1,T[$, first order conditions read as
\[
\begin{array}{ll}
\displaystyle
0=-\frac{1}{2}\left(\frac{\theta_S-\theta_1}{\overline\tau-\underline\sigma_1}\right)^2+{\cal F}^{1,1}({\cal M}(\overline\tau))+V'(\theta_1,\overline\tau,1,1)=\\
\displaystyle
-\frac{1}{2}\left(\frac{\theta_S-\theta_1}{\overline\tau-\underline\sigma_2}\right)^2+{\cal F}^{1,1}({\cal M}(\overline\tau))+V'(\theta_1,\overline\tau,1,1)
\end{array}
\]
contradicting   $\underline\sigma_1\neq \underline\sigma_2$. Note that $V$ is not necessarily differentiable in time, however it has a super-differential at any instant. Indeed it is easy to see that the value function for $(w_1,w_2)=(0,0)$ has such a property, and then obtain the super-differentiability of the others arguing backward in (\ref{eq:V00})(\ref{eq:V01}) (\ref{eq:V10}) (\ref{eq:V11}).
The super-differentiability at $\overline\tau$ implies that, at least locally in time around $\overline\tau <T$, one
$V(\theta_1,\tau,1,1)\le h(\tau)$ where $h$ is a suitable differentiable function, with equality holding at $\overline\tau$. Hence the argument above would hold with $V$ replaced by $h$.

(2) Let us  suppose that agents arriving at  time $\underline \sigma$ in  $(\theta_S,1,1)$ and moving towards $\theta_1$ may reach this latter significant state in more than one optimal time, say $\overline \tau_1,\overline \tau_2$  with $\overline \tau_1< \overline \tau_2$. 
This means that 	
$\underline \sigma(\theta_1,\overline\tau_1)=\underline \sigma(\theta_1,\overline\tau_2)=\underline \sigma$.
We now observe that only agents entering $\theta_S$ at time $\underline \sigma$ may reach $\theta_1$ between $\tau_1$ and $\tau_2$,  as  agents cannot overtake each other (easy to prove). As a consequence $\underline \sigma(\theta_1,s)$ is constant in the interval $[\tau_1,\tau_2]$ and hence its time derivative $\underline \sigma_s(\theta_1,s)$ is null. This last fact  in turn implies that the value of the entrance flow at $\underline \sigma$, i.e., $g(\underline \sigma)$, would  uniformly spread over the interval $[\tau_1,\tau_2]$ and hence would become equal to 0.

Remark that the above arguments imply that no Dirac masses can arise at any point and at any time, apart from the case of a significant point at the final time $T$, or, just after a switching, when the new choice of the optimal control is $u=0$, i.e. to not move. Both situations do not affect the flow functions $g_{i,j}$.

\section{Appendix: on HBJ and transport equations}
For the reader convenience, we here recall the 
Hamilton-Jacob-Bellman (3)--(6) and the transport equations (8a)--(8d) contained in \cite{bagpes}.  {We point out once  again that the final cost in the present paper is different from that in \cite{bagpes}: here  $c_S\xi_{\theta=\theta_S}(T)$ replaces $c_3(\theta-\theta_S)^2$.  

In \citet[Theorem 1]{bagpes} was proved that the value function $V(\theta,t,w_1,w_2)$ %\deleted{described by  (\ref{eq:V00}), (\ref{eq:V01}), (\ref{eq:V10}), (\ref{eq:V11}),} 
 is a continuous and bounded viscosity solution of the following   HJB equations (subscript indicate derivatives).

\noindent In $B_{1,1}$:
\begin{subequations}\label{eq:HBJXX}
	\begin{align}
	\label{eq:HJB11}
	&	\left\{
	\begin{array}{ll}
	\displaystyle
	-V_t(\theta,t,1,1)+\frac{1}{2}|V_\theta(\theta,t,1,1)|^2={\cal F}^{(1,1)}({\cal M}(t))&\mbox{in } ]\theta_1,\theta_2[\times]0,T[\\
	\displaystyle
	V(\theta_1,t,1,1)=V(\theta_1,t,0,1)&\mbox{in } ]0,T]\\
	\displaystyle
	V(\theta_2,t,1,1)=V(\theta_2,t,1,0)&\mbox{in } ]0,T]\\
	\displaystyle
	V(\theta,T,1,1)=c_1+c_2+c_3(\theta-\theta_S)^2&\mbox{in } ]\theta_1,\theta_2[
	\end{array}
	\right.\\
	\intertext{In   $B_{0,1}$:}
	\label{eq:HJB01}
	&	\left\{
	\begin{array}{ll}
	\displaystyle
	-V_t(\theta,t,0,1)+\frac{1}{2}|V_\theta(\theta,t,0,1)|^2={\cal F}^{(0,1)}({\cal M}(t))&\mbox{in } ]\theta_2-2\pi,\theta_2[\times]0,T[\\
	\displaystyle
	V(\theta_2-2\pi,t,0,1)=V(\theta_2,t,0,0)&\mbox{in } ]0,T]\\
	\displaystyle
	V(\theta_2,t,0,1)=V(\theta_2,t,0,0)&\mbox{in } ]0,T]\\
	\displaystyle
	V(\theta,T,0,1)=c_2+c_3(\theta-\theta_S)^2&\mbox{in } ]\theta_2-2\pi,\theta_2[
	\end{array}
	\right.\\
	\intertext{In   $B_{1,0}$:}
	\label{eq:HJB10}
	&	\left\{
	\begin{array}{ll}
	\displaystyle
	-V_t(\theta,t,1,0)+\frac{1}{2}|V_\theta(\theta,t,1,0)|^2={\cal F}^{(1,0)}({\cal M}(t))&\mbox{in } ]\theta_1,\theta_1+2\pi[\times]0,T[\\
	\displaystyle
	V(\theta_1,t,1,0)=V(\theta_1,t,0,0)&\mbox{in } ]0,T]\\
	\displaystyle
	V(\theta_1+2\pi,t,1,0)=V(\theta_1,t,0,0)&\mbox{in } ]0,T]\\
	\displaystyle
	V(\theta,T,1,0)=c_1+c_3(\theta-\theta_S)^2&\mbox{in } ]\theta_1,\theta_1+2\pi[
	\end{array}
	\right.\\
	\intertext{and in  $B_{0,0}$:}
	\label{eq:HJB00}
	&	\left\{
	\begin{array}{ll}
	\displaystyle
	-V_t(\theta,t,0,0)+\frac{1}{2}|V_\theta(\theta,t,0,0)|^2={\cal F}^{(0,0)}({\cal M}(t))&\mbox{in } \mathbb{R}\times]0,T[\\
	\displaystyle
	V(\theta,T,0,0)=c_3(\theta-\theta_S)^2&\mbox{in } [0,2\pi].
	\end{array}
	\right.
	\end{align}
\end{subequations}

\bigskip

Moreover, the four transport equations for the density $m$, one per every branch, are
\begin{subequations}\label{eq:transportXX}
	\begin{align}
	\label{eq:transport11}
	&	\left\{
	\begin{array}{ll}
	\displaystyle
	m^{1,1}_t(\theta,t)-[V_\theta(\theta,t,1,1)m^{1,1}(\theta,t)]_\theta=0\ \mbox{in } B_{1,1}\times[0,T]\\
	\displaystyle
	m^{1,1}(\theta_S,t)=g(t)
	\end{array}
	\right.\\
	\nonumber\\
	\label{eq:transport10}
	&	\left\{
	\begin{array}{ll}
	\displaystyle
	m^{1,0}_t(\theta,t)-[V_\theta(\theta,t,1,0)m^{1,0}(\theta,t)]_\theta=0\ \mbox{in } B_{1,0}\times[0,T]\\
	\displaystyle
	m^{1,0}(\theta_2,t)=m^{1,1}(\theta_2,t)
	\end{array}
	\right.\\
	\nonumber\\
	\label{eq:transport01}
	&	\left\{
	\begin{array}{ll}
	\displaystyle
	m^{0,1}_t(\theta,t)-[V_\theta(\theta,t,0,1)m^{0,1}(\theta,t)]_\theta=0\ \mbox{in } B_{0,1}\times[0,T]\\
	\displaystyle
	m^{0,1}(\theta_1,t)=m^{1,1}(\theta_1,t)
	\end{array}
	\right.\\
	\nonumber\\
	\label{eq:transport00}
	&	\left\{
	\begin{array}{ll}
	\displaystyle
	m^{0,0}_t(\theta,t)-[V_\theta(\theta,t,0,0)m^{0,0}(\theta,t)]_\theta=0\ \mbox{in } B_{0,0}\times[0,T]\\
	\displaystyle
	m^{0,0}(\theta_1,t)=m^{1,0}(\theta_1,t)+m^{1,0}(\theta_1+2\pi,t)\\
	\displaystyle
	m^{0,0}(\theta_2,t)=m^{0,1}(\theta_2,t)+m^{0,1}(\theta_2-2\pi,t).
	\end{array}
	\right.
	\end{align}
\end{subequations}
We recall that, consistently with  \eqref{eq:optimalfeedback}, one has  
$u^*=-V_\theta$.

\end{document}